\documentclass[11pt,a4paper]{article}

\usepackage[a4paper,top=1.2in, bottom=1.2in, left=1in, right=1in]{geometry}
\usepackage{commath}
\usepackage{amsmath,paralist,amsthm,mathpazo}
\usepackage{amsfonts,amssymb}
\usepackage{lmodern}
\usepackage{ulem}
\usepackage{setspace}
\usepackage{enumerate}
\usepackage{pgfplots}
\usepackage{float}
\usepackage{enumitem}
\usepackage{tikz}
\usepackage{multirow,array}
\usepackage{pbox}
\usepackage{graphicx}
\usepackage[utf8]{inputenc}
\usepackage{amssymb}
\usepackage{adjustbox}
\usepackage{fancyref}
\usepackage[utf8]{inputenc}
\usepackage{scalerel}
\usepackage[round]{natbib}
\usepackage{lineno}
\usepackage{setspace}
\usepackage{subfigure}
\usepackage{tikz}
\usetikzlibrary{positioning}
\usepackage{float}
\usepackage{caption}
\usepackage{tabstackengine}
\usepackage{graphics,graphicx}
\usepackage{pst-node,pst-tree,pst-pdf,pdftricks,pstricks-add}
\usepackage{pstricks}

\graphicspath{ {E:/images/} }

\usepackage[justification=centering]{caption}
\usetikzlibrary{positioning}

\newtheorem{proposition}{Proposition}

\newtheorem*{remark}{\textit{Remark}}

\newtheorem{conjecture}{Conjecture}

\newtheoremstyle{as}
{\topsep}
{\topsep}
{\itshape}
{}
{\scshape}
{.}
{ }
{\thmname{\textbf{\textit{#1}}}\thmnumber{ $\mathbf{(#2)}$}\thmnote{ (#3)}}       
\theoremstyle{as}

\newtheoremstyle{pro}
{\topsep}
{\topsep}
{\itshape}
{}
{\scshape}
{.}
{ }
{\thmname{\textbf{\textit{#1}}}\thmnumber{ $\mathbf{#2}$}\thmnote{ (#3)}}%
\theoremstyle{pro}

\newtheoremstyle{lem}
{\topsep}
{\topsep}
{\itshape}
{}
{\scshape}
{.}
{ }
{\thmname{\textbf{\textit{#1}}}\thmnumber{ $\mathbf{#2}$}\thmnote{ (#3)}}%
\theoremstyle{lem}

\usepackage{subfigure}
\usepackage{tikz}

\begin{document}
\centerline{\textbf{\Large The Uniformed Patroller Game}}
\bigskip
\centerline{\textbf{Steve Alpern}\footnote{\scriptsize ORMS Group, Warwick Business School, University of Warwick, Coventry CV4 7AL, UK, steve.alpern@wbs.ac.uk},\textbf{Stamatios Katsikas}\footnote{\scriptsize Centre for Complexity Science, University of Warwick, Coventry, CV4 7AL, UK,  stamkatsikas@gmail.com}}
\smallskip
\centerline{October $13$, $2018$}
\bigskip
\begin{abstract}
\noindent  
In the recently introduced \textit{network patrolling game}, an \textit{Attacker} carries out an attack on a node of her choice, for a given number $m$ of consecutive periods. The parameter $m$ indicates the difficulty of the attack at a given node. To thwart such an attack, the \textit{Patroller} adopts a walk on the network, hoping to be at the attacked node during one of the $m$ periods. If this occurs, the attack is intercepted and the Patroller wins; otherwise the Attacker wins. In the original setting, the Attacker has no knowledge of the Patroller's location at any time. Here, to model the important alternative where the Patroller can be identified when he is at the Attacker's node (e.g. the Patroller wears a uniform), we allow the Attacker to initiate her attack after waiting for a chosen number $d$ of consecutive periods in which the Patroller has been away. We solve this version of the game for various networks: star, line, circle and a mixture. We restrict the Patroller to Markovian strategies, which cover the whole network. 
\end{abstract}

\section{Introduction}\label{Section1}

When patrolling a network against an attack, or infiltration, at an unknown node, there are two plausible scenaria. Either $(i)$ the Patroller may be essentially invisible to the Attacker (e.g. a stealthy drone, a plainclothes policeman, or a driver of an unmarked car); or $(ii)$ he may be identifiable by the Attacker (e.g. through his uniform - hence our title, or his driving a police car). Thus far the literature on patrolling games has concentrated exclusively on the first (invisible) scenario. This is the first attempt to model the second scenario of a uniformed Patroller. Some consequences of this distinction are fairly obvious; for example, a rational mugger will never start an attack in a subway car when a uniformed policemen is present. Other consequences, to be explored here, are less obvious. 

A game theoretic model of the Attacker-Patroller conflict on a network has recently been the subject of several investigations. \cite{Alp-1} modeled the problem as a zero-sum game between a Patroller and an Attacker, where the Attacker can attack a chosen node in a chosen time period, and the Patroller hopes to intercept the attack, by following a chosen walk on the network. In this zero-sum game the payoff to the maximizing Patroller is the probability that the attack is intercepted. Formally, the Attacker chooses a node $i$ to attack for a time interval $J$, that is a sequence of $m\geq 2$ consecutive time periods, and the attack is intercepted if the Patroller is at node $i$ at some time period within the interval $J$. Here, we keep these dynamics and payoff the same. However, to simplify the problem we restrict the Patroller to ergodic Markovian strategies, where every time he arrives at a node he leaves it by the same distribution over its neighbors.

Note that in the original formulation given above, the Attacker's strategy $(i,J)$ cannot depend in any way on the Patroller's locations at any time prior to $J$. However, in many real world situations of this type, the Patroller is identifiable; he might be wearing a uniform, or driving in a marked car. Thus, as we assume in our model presented here, the Attacker can go to the node she wishes to attack, and wait there prior to her attack. She can observe when the Patroller is present at her chosen node, and when he is not. Once the Patroller leaves her node, the Attacker only knows that the Patroller is not there; she cannot see the Patroller from a distance (vision is limited to the attack node). This assumption defines a game $U(Q,m)$, where $Q$ is the network to be defended, and $m$ is the difficulty of attack or infiltration as measured in terms of the time required. 

The Attacker's pure strategy can be defined as a pair $(i,d)$, where she goes to node $i$ and initiates her attack after the Patroller has been there first and then has been away for $d$ consecutive periods. Since the Patroller wears a uniform, his Markovian strategy is assumed to be known to the Attacker, so we adopt a Stackelberg approach where the Patroller is the first to move. This is a common assumption in patrolling problems. We solve this game for small values of $m$ and several families of networks: star, line, circle, and star-in-circle. In particular, both line and circle networks can be interpreted as perimeters of regions to be defended against infiltration. As such, these are models of border defense by uniformed patrollers, modifying the border patrol game of \cite{Pap-1}. 

The motivation for \textit{the Patroller in Uniform Problem} came to the first author a long time ago, when in the 1970's a uniformed policeman was assigned to every subway train (consisting of ten cars, thus being equivalent to the line network $L_{10})$ in New York City. In the first months at least, the policeman patrolled in a back and forth motion, from car $\mathit{1}$ to car $\mathit{10}$ and backwards. This patrol seemed quite foolish as, in our notation, the attack strategy $(\mathit{1},1)$ (attacking an end car as soon as the Patroller has left it) would guarantee a win even if the difficulty $m$ of the attack (i.e. the attack duration) was as large as $17$. Finding optimal patrolling strategies in this context remains an open issue until now, over thirty years later. Note that if the Patroller on this train is following a random walk and the time required for a mugging is say $m=2,$ it is not optimal to attack as soon as the Patroller leaves your end car (i.e. with $d=1$), as he will catch you with probability $1/2$. It is clearly better to wait for some larger number of consecutive periods $d$ in which the Patroller is away.

Patrolling problems have been studied for long, see e.g. \cite{Morse}, but only from the Patroller's viewpoint. A game theoretic approach, modelling an adversarial Attacker who wants to infiltrate or attack a network at a node of her choice, has only recently been introduced by \cite{Alp-1}. The techniques developed there were applied to the class of line networks by \cite{Pap-1}, with the interpretation of patrolling a border. Other research following a similar reasoning includes \cite{Lin-1} for random attack times, \cite{Lin-2} for imperfect detection, \cite{Hochbaum} on security routing games, and \cite{Basilico-2} for uncertain alarm signals. See also \cite{Baykal} on infracstructure security games. Earlier work on patrolling a channel/border with different paradigms, includes \cite{Washburn-1,Washburn-2}, \cite{Szechtman}, \cite{Zoroa}, and \cite{Collins}. The case where the Patroller is restricted to periodic walks has recently been studied by \cite{Alp-2}. The related problem of ambush is studied by \cite{Baston-1,Baston-2}, while an artificial intelligence approach to patrolling is given by \cite{Basilico-1}. Applications to airport security and counter terrorism, are given respectively by \cite{Pita}, and \cite{Fokkink}. The problem studied here, of patrolling when the Attacker can see the Patroller when he is in close proximity, appears to be new.

Our approach is also related to two problems studied in the theoretical computer science literature. The first one is the \textit{cyclic routing problem} of \cite{Ho} which, in the case of a single patroller, asks if the patroller can return to every node $i$ within at most $\kappa_{i}$ periods. This is a variant of the \textit{pinwheel problem} of \cite{Holte}, and the \textit{patrolling security games} as studied in \cite{Basilico-1} and \cite{Basilico-2}, and is akin to the problem faced by a waiter in a restaurant who must return to each table with a certain frequency depending on the number of diners at the table. Another version of patrolling games where the two players can see each other at all times (really a pursuit-evasion game) is the \textit{spy game}, as in \cite{Cohen}.

The paper is organized as follows. In Section \ref{Section2} we explain our model giving a specification of the Patroller in uniform problem on an arbitrary finite network. In Section \ref{Section3} we solve the game for the star network $S_{n}$ with a central node connected to $n$ ends. In this setting we also examine an extension of our initial framework by partially loosening the Markovian restriction on Patroller's movement. In Sections \ref{Section4} and \ref{Section5} we solve the game for the line network $L_{n}$ and the circle network $C_{n}$ respectively, for $n=4,5$, using numerical methods. In Section \ref{Section6} we consider a hybrid network consisting of a circle network with a center connected to all ends. In the last setting we examine another extension of our initial framework by allowing the Attacker to observe towards which direction the Patroller leaves the attack node. In Section \ref{Section7} we summarize our conclusions and suggest areas of future research.

\section{Formal Model}\label{Section2}

An attack strategy is a pair $(i,d)$,  where $i$ is a node of a given finite network $Q$, and $d$ is the number of consecutive time periods that the Patroller must be away from node $i$ for the Attacker to initiate her attack (we use for this alternately the terms \textit{waiting time} and \textit{delay}). In \cite{Alp-1} there is a time horizon denoted by $T$ which denotes the last period that the attack has to finish by. By analogy, here we introduce a positive integer $D$ which is an upper bound for the delay $d$, that is, $d\leq D$. The value of $D$ is common knowledge. A strategy for the Patroller is an ergodic Markov chain on the nodes of the given network. We require that this strategy respects all symmetries of the network. For example on a circle network the Patroller must move clockwise and counter clockwise with equal probability. 

To illustrate the syntax of the Attacker's actions, and in particular the waiting parameter $d$, consider the following example. Suppose that the Patroller's presence at the attack node $i$ is indicated by a $1$, respectively his absence is indicated by a $0$. Suppose further that the Attacker waits until the Patroller arrives there, and after that his presence-absence sequence that she records is $11010100\dots$. If the Attacker's waiting time is for example $d=2$, then her decision as to whether or not to attack in each subsequent period is illustrated in Table \ref{Table1}.

\setlength{\tabcolsep}{13pt}
\renewcommand{\arraystretch}{1}
\begin{table}[h!]
\centering
\begin{tabular}{ c  c  c  c  c  c  c  c }
\hline\hline
$1$ & $11$ & $110$ & $1101$ & $11010$ & $110101$ & $1101010$ & $11010100$ \\ \hline
\textit{wait} & \textit{wait} & \textit{wait} & \textit{wait} & \textit{wait} & \textit{wait} & \textit{wait} & \textit{attack!} \\ 
\hline\hline
\end{tabular}
\caption{The Attacker attacks after $d=2$ periods of the Patroller's absence}
\label{Table1}
\end{table}

If for example the attack difficulty is $m=3$, then the attack will be successful only if the Patroller's sequence continues with two more $0$'s, that is, $1101010000\dots$. Note that the attack can begin in the same period that the Patroller's absence has been observed. Thus, we consider $m\geq 2$ as otherwise the Attacker could win simply by attacking as soon as the Patroller is not present at the attack node. 

\begin{figure}[H]
\centering
\includegraphics[scale=0.85]{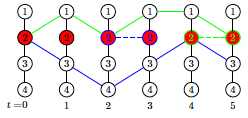}
\caption{Example of Attacker-Patroller dynamics for the line graph $L_{4}$.}
\label{Figure1}
\end{figure}

The above discussion explains the waiting parameter $d$ from the Attacker's point of view, who simply observes the Patroller's presence or absence at her node. To illustrate the attack mechanism in the context of both players, take the network $Q$ as the line graph $L_{4}$ with four nodes, labeled $1$ to $4$ from top to bottom, and consider the dynamics presented in Figure \ref{Figure1}, where the time axis for the first five periods is drawn horizontally, and the line graph is drawn vertically.

Suppose that the Attacker chooses to attack node $\mathit{2}$ (drawn in red), of difficulty $m=2$, with delay $d=2$, and the Patroller is there at time $t=0$. We consider two patrols adopted by the Patroller on the line network $L_{4}$; the walk $w_{1}=(2,1,2,1,1,2,\dots)$ drawn in green on top, and the walk $w_{2}=(2,3,4,3,2,3,\dots)$ drawn in blue on the bottom. Take first the Attacker's response to $w_{1}$. At time $t=2$ she resets her waiting clock to zero, but at time $t=4$ since the Patroller is still away after $d=2$ periods, she initiates her attack that lasts for the interval $\{4,5\}$ indicated by a dashed green horizontal line. However, since the Patroller is back at node $\mathit{2}$ at time $t=5$, the attack is intercepted (the green patrol intersects the dashed green horizontal attack). Next, consider the Attacker's response to $w_{2}$. At time $t=2$ she begins the attack that lasts for the interval $\{2,3\}$ indicated by a dashed blue horizontal line. Since the Patroller is not back at node $\mathit{2}$ after $m=2$ periods, the attack is not intercepted (the blue patrol is disjoint from the dashed blue horizontal attack). Thus, against the attack strategy (node $ \mathit{2}$, $d=2)$, the patrol $w_{1}$ wins for the Patroller, while the patrol $w_{2}$ loses.

The reader might as well think that the Attacker has additional strategies which are not specified by our restriction to pairs $(i,d)$. For example the Attacker could attack after say three consecutive $0$'s when she initially arrives at node $i$, without necessarily first waiting for the Patroller to visit node $i$ and then counting three consecutive $0$'s. However we show that the Attacker can always gain at least as good an outcome using some strategy $(i,d)$. In particular, since we assume that the Patroller has been patrolling for an arbitrarily long time before the period when the Attacker arrives at node $i$, then we could tell the Attacker, for free, the total number of $0$'s at node $i$ since the last $1$ (including the three $0$'s the Attacker has witnessed). Say for example this total is $7$. If the Attacker chooses to ignore this new information, and attack as initially planned, then her expected payoff is the same as that of the strategy $(i,7)$. On the other hand, if the Attacker changes her strategy and decides not to attack, then she is again consistent with some strategy $(i,d)$.

We will be mainly concerned with specific networks in this paper. However, we make the following conjecture regarding general networks.

\begin{conjecture}\label{Conjecture1}
For any network and any attack duration, it is optimal to always reflect at leaf~nodes.
\end{conjecture}

Note that we establish Conjecture \ref{Conjecture1} in certain cases in the paper (see Proposition \ref{Proposition2}).

\subsection{Results from Patroller without a uniform}\label{Section2.1}
\indent

For the most part, the techniques required to solve the patrolling game with uniform are distinct from those developed in \cite{Alp-1} without the uniform. However, in a few cases reviewed here, some of those results carry over, possibly with restrictions, to the case of a uniform. In the context of non uniformed patrolling games the Patroller chooses any walk on the graph, while the Attacker chooses the time interval and node for the attack.

\begin{figure}[h!]
 \begin{minipage}[b]{0.45\textwidth}
    \centering
\begin{figure}[H]
\centering
\includegraphics[scale=0.85]{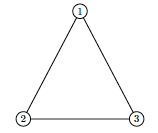}
\caption{The Circle network $C_{3}$.}
\label{Figure2}
\end{figure}
  \end{minipage}
\hfill
  \begin{minipage}[b]{0.49\textwidth}
    \centering
\begin{figure}[H]
\centering
\includegraphics[scale=0.85]{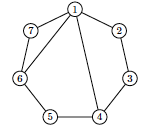}
\caption{An asymmetric Hamiltonian graph.}
\label{Figure3}
\end{figure}
    \end{minipage}
\end{figure}

The elementary Theorem $13$ of the earlier paper observes that when the graph $Q$ is Hamiltonian with $n$ nodes, the value of the game is given by $m/n$ (or by $1$, if $m\geq n$), where $m$ is  the attack duration. This is achieved by the Patroller following a Hamiltonian cycle endlessly, starting at a random time mod $n$. Such a Patroller strategy can be implemented in a Markovian manner by choosing the out edge in the Hamiltonian path with probability $1$. However, when $m\leq n-1$, the Attacker in the uniformed context can counter this strategy by attacking as soon as the Patroller leaves her node (i.e. the attack node), and completing her attack before the Patroller returns. When $m=n$ the Attacker has no counter move to this strategy, and the value is $1$ (as long as the Hamiltonian path strategy does not violate our requirement that the Patroller respects the symmetry of the graph). For the graph in Figure \ref{Figure3} ($n=7$), the interception probability when the Patroller goes around the circle say clockwise is equal to $1$ when $m=7$. (This graph has no symmetries, that is, the isomorphism group is the identity alone. To see this note that nodes cannot go into nodes of different degrees, so that $1\longleftrightarrow 1$ and hence any non identity isomorphism satisfies $6\longleftrightarrow 3$, which violates degree preservation). For the Hamiltonian graph in Figure \ref{Figure2}, where $C_{3}$ is a cycle graph and we take $m=3$, symmetry prevents the Patroller from distinguishing between the two Hamiltonian cycles. Actually this is the complete graph $K_{3}$ on three nodes, for which, as we show in equation \eqref{eq34} of Section \ref{Section6.1}, the value is equal to $3/4$, with the random walk optimal for the Patroller and an arbitrary delay $d$ optimal for the Attacker.

Prompted by this problem, we note an extension of Theorem $13$ of \cite{Alp-1} which has the same proof as given above: \textit{Consider the (non uniformed) patrolling game on a graph $Q$ for which the length $\kappa$ of a minimal covering cycle satisfies $\kappa\leq m$. Then the value (i.e. the optimal interception probability) is $\mathit{1}$.} This in fact is a useful result for the original patrolling game.

The assumption in the uniformed patrolling game that the Patroller's strategies respect all symmetries in the graph is not essential. A more general approach would give, as a parameter of the~game, not only the graph $Q$, but also some symmetry group of $Q$ (possibly the identity alone) that had to be respected. This approach was adopted for the rendezvous searchers in \cite{Alp}. For example, for one of these symmetry groups two distinguishable rendezvous searchers on a continuous circle can optimize by going clockwise and anticlockwise. However, given other symmetry groups they cannot have a common understanding of clockwise. These notions are explained at length in that paper.

Of course, if the Patroller can remember how he got to a particular node (one step memory), this is enough to propel him around the cycle in a consistent direction, as noted on page $1256$ of \cite{Alp-1}: \textit{``An absent-minded Patroller who randomizes anew at each time unit depending on her location could be modeled through a Markov game. Such a game may have interesting features: for example, a Patroller who can remember from which direction she has come may perform better on a cycle rather than a complete graph, because she is less likely to backtrack, contradicting the ``increasing in edges'' property''}.

Another result which cannot be transferred to the case of a Patroller in uniform is Lemma $5$ of \cite{Alp-1}, which says that attacks on penultimate nodes (i.e. nodes adjacent to a leaf node) are dominated. For example, even though we find numerically for line graphs that attacks at the endpoints (the so called ``diametrical strategy'' of Lemma $9$ of \cite{Alp-1}) are optimal, it seems that general results such as those quoted as Lemmas $5$ and $9$ above are much more difficult to obtain for a Patroller in a uniform rather than without it.

\section{The Star Network $S_{n}$}\label{Section3}

The network $S_{n}$ is the star with a single center connected to $n$ end nodes. We restrict the Patroller to Markovian strategies that reflect from the ends with probability $s$, from the center $c$ go to each end with equal probability $p$, and remain at the center with probability $r=1-n\cdot p$. This setting simplifies the game by introducing a single parameter family of patrolling strategies.

\begin{figure}[h!]
\centering
\includegraphics[scale=0.85]{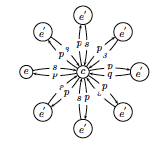}
\caption{The Star network $S_{n}$.}
\label{Figure4}
\end{figure}

We begin by assuming that the attack takes place at an end node, which we denote by $e$, and then we show that the Patroller should reflect from the ends ($s=1$). Note that since we have taken $m\geq 2$, reflecting from the ends will further imply that the Attacker should never attack at the center $c$ because in that case the Patroller will never be away from it for two consecutive periods.

\subsection{Attack duration $m=2$}\label{Section3.1}
\indent

Initially we assume that an attack takes $m=2$ periods. In this case, an attack at node $e$ cannot be intercepted if it starts when the Patroller is at another end node $e^{'}\neq e$, but only if the Patroller is at the center. Therefore, an attack at node $e$ will be intercepted with probability $q\cdot p$, where $q$ is the probability that the Patroller is at the center $c$ at the beginning of the attack, which in turn implies that the Attacker should choose the waiting time $d$ so as to minimize $q$.

We wish to calculate how the probability $q$ that the Patroller is at $c$ changes over time, as the Patroller continues to be away from the Attacker's chosen node $e$. So suppose that at some period $t$ the Patroller is not at $e$, but he is either at $c$ with probability $q$, or at one of the end nodes other than $e$ with probability $1-q.$ Then, in the following period $t+1$ the Patroller will be either at $c$ with probability $q\cdot r +\left( 1-q\right)\cdot s$, or he will be at node $e$ with probability $q\cdot p +(1-q)\cdot 0$. 

Hence, conditional on the Patroller not being at node $e$ at period $t$, the probability that he is at the center $c$ at period $t+1$ is given by

\begin{equation}\label{eq1}
f(q,s)=\frac{q\cdot r+(1-q)\cdot s}{1-p\cdot q}.
\end{equation}

Fraction \eqref{eq1} is increasing in $s$, and since it is $s\leq 1$, then we have that $f(q,s)$ is maximized for $s=1$ at

\begin{equation}\label{eq2}
f(q)=f(q,1)=\frac{q\cdot r+(1-q)}{1-p\cdot q}=\frac{1-n\cdot p\cdot q}{1-p\cdot q}.
\end{equation}

Additionally, since from equation \eqref{eq2} we find that $f^{'}(q)=-\frac{p\cdot (n-1)}{(1-p\cdot q)^{2}}<0$, then $f(q)$ is decreasing and hence minimized for $q=1$, with $f(1)=\hat{q}$, where

\begin{equation}\label{eq3}
\hat{q}=\frac{(1-n\cdot p)}{(1-p)}.
\end{equation}

The Attacker can obtain this minimum probability $\hat{q}$ of the Patroller being at node $c$ at the beginning of her attack, by initiating her attack on the second period that the Patroller is away from her chosen node $e$, that is, by adopting the waiting time $d=2$. Notice that the optimal Attacker strategy $(e,2)$ does not depend on $p$, so the assumption that the Attacker knows $p$ is not necessary in this case. The attack $(e,2)$ will be intercepted if the Patroller is at $c$ in its first period and he goes to $e$ in its second period, that is, with probability

\begin{equation}\label{eq4}
a_{2}(n,p)=\hat{q}\cdot p=b(p)=\frac{(1-n\cdot p)\cdot p}{1-p}.
\end{equation}

\begin{figure}[h!]
	\centering
\includegraphics[scale=0.425]{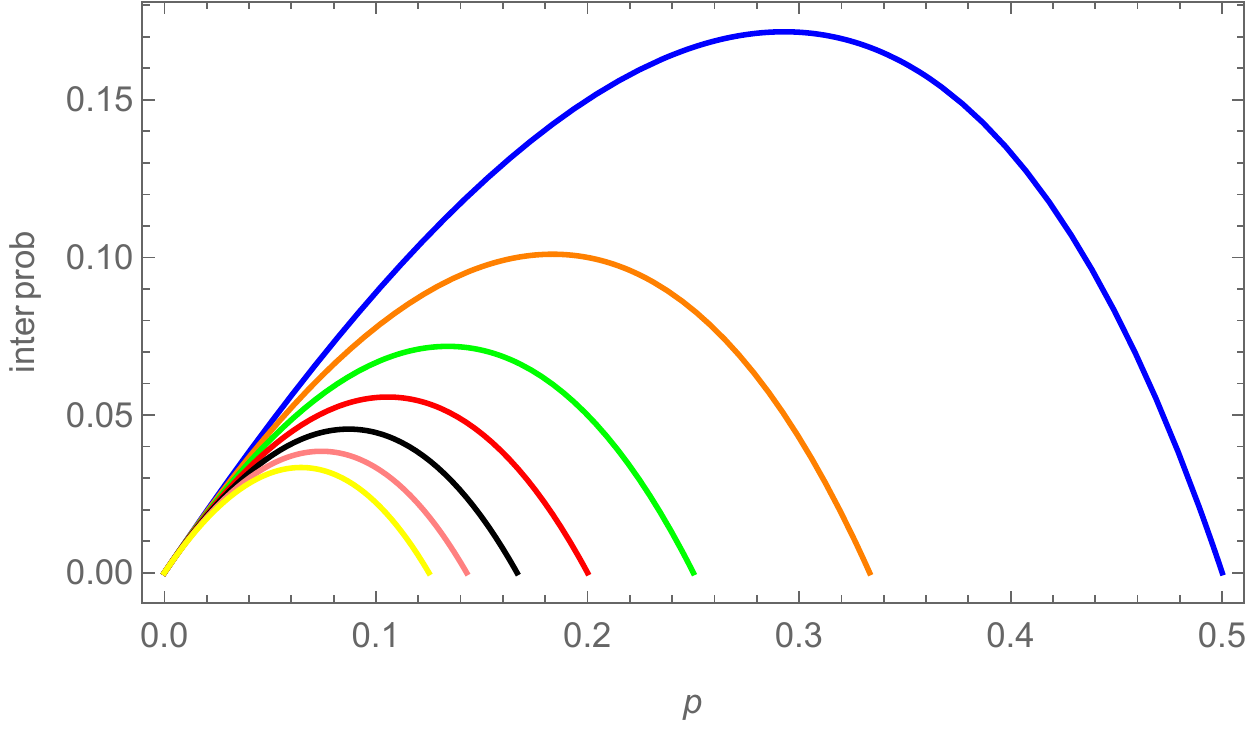}
\caption{The interception probability $a_{2}(n,p)$, for $n=2(\text{blue}), \dots, 8(\text{yellow})$.}
\label{Figure5}
\end{figure}

For a given star network $S_{n}$, the Patroller will choose the value of $p$ that maximizes the interception probability \eqref{eq4}. Figure \ref{Figure5} shows the variation of \eqref{eq4} with $p$, respectively for $n=2,\dots,8$ arcs in $S_{n}$, where $p\in [0,1/n]$ is the probability with which the Patroller moves to an end from the center $c$.

For the $n=2$ arc star (which is equivalent to a line network with three nodes), the Patroller should move from the center towards each end with probability about $0.3$, and remain at the center with probability about $0.4$. For the $n=3$ arc star, the Patroller should move towards each of the three end nodes with probability about $0.2$, and remain at the center with probability about $0.4$, etc.. 

We can be more precise about these Patroller's optimal strategies. In particular, the optimal value $\hat{p}=\hat{p}(n)$ for $p$ depends on $n$ and can be found by solving the first order equation

\begin{equation*}
a_{2}^{'}(p)=\frac{(n\cdot p^{2}-2\cdot n\cdot p+1)}{(1-p)^{2}}.
\end{equation*}

This can be simply written in the form

\begin{equation*}
n\cdot p^{2}-2\cdot n\cdot p +1=0 ,
\end{equation*}

\noindent giving

\begin{minipage}{0.47\linewidth}
\begin{equation}\label{eq5}
\hat{p}=1-\frac{\sqrt{n\cdot (n-1)}}{n}
\end{equation}
\end{minipage},
\begin{minipage}{0.47\linewidth}
\begin{equation}\label{eq6}
\hat{r}=\sqrt{n\cdot (n-1)}-(n-1)
\end{equation}
\end{minipage}.

Equations \eqref{eq5}, \eqref{eq6} show respectively that the optimal probability $\hat{p}$ is asymptotic to $1/(2\cdot n)$, while the optimal probability $\hat{r}$ goes to $1/2$. The value $V$ of the game is given by

\begin{equation}\label{eq7}
V=a(n,\hat{p})=(2\cdot n -1)-2\cdot\sqrt{n\cdot (n-1)} .
\end{equation}

\begin{figure}[h!]
	\centering
	\subfigure[Patroller's optimal walk $\hat{r}(n)$]{\includegraphics[scale=0.425]{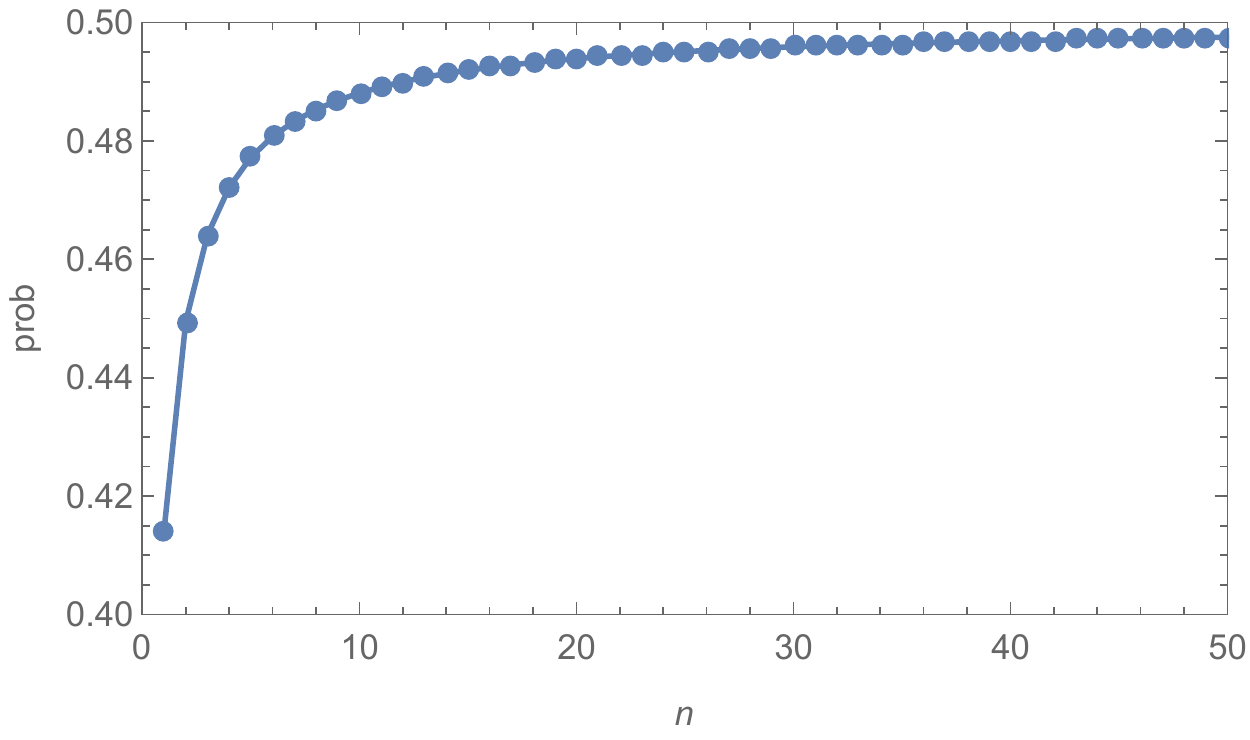}}
\hspace{1.75em}
           \subfigure[interception probability $V(n)$]{\includegraphics[scale=0.425]{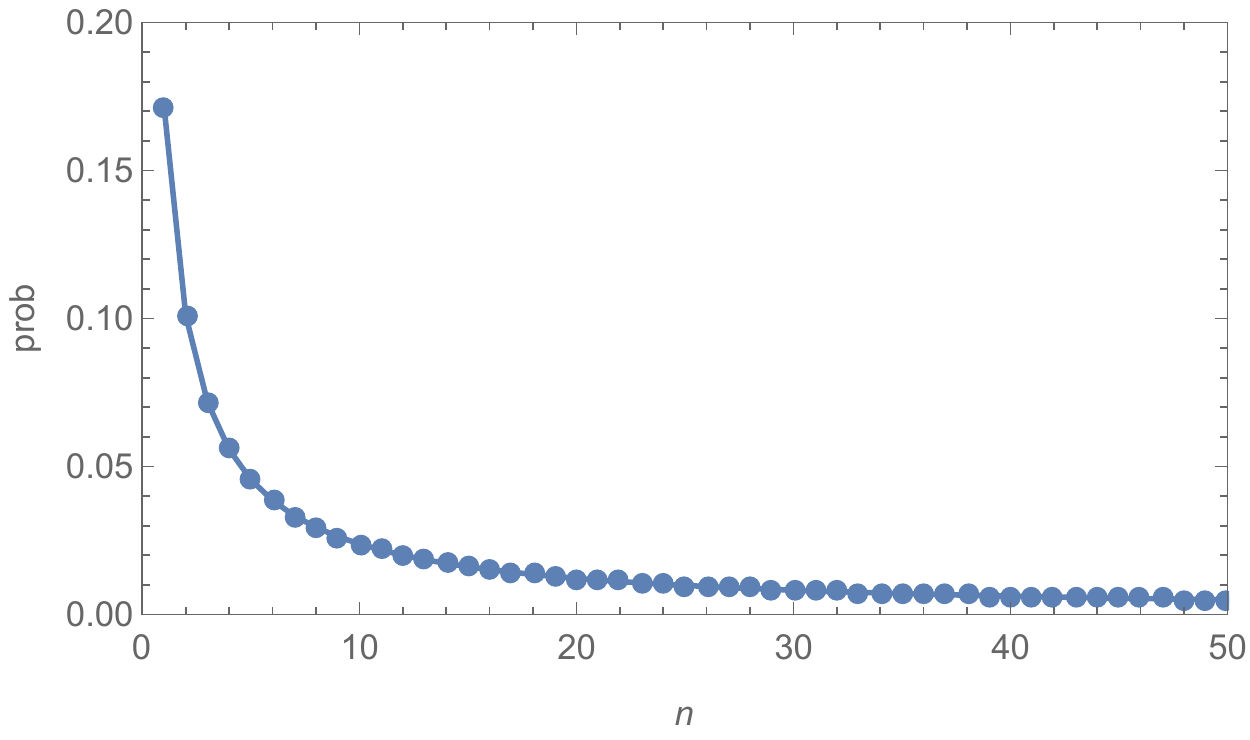}}
	\caption{Optimal values for $m=2$.}
	\label{Figure6}
\end{figure}

In Figure \ref{Figure6} we plot the optimal Patroller's strategy, in terms of the probability $\hat{r}$ of remaining at the center, and the corresponding value of the game for increasing number of nodes $n$.

We can pull together the above results into the following proposition:

\begin{proposition}\label{Proposition1}
Consider the Uniformed Patroller Problem on the star network $S_{n}$ with $n$ arcs and an attack difficulty of $m=2$. The optimal Patroller's strategy is to reflect off the end nodes and to remain at the central node with probability given by \eqref{eq6}. The optimal Attacker's strategy is to locate at a random end node and  begin the attack in the second period that the Patroller is away $(d=2)$. The interception probability under these strategies is given by \eqref{eq7}.
\end{proposition}

\begin{remark}\label{Remark1}
Note that Proposition \ref{Proposition1} gives an attack strategy which does not require any knowledge of the Patroller's Ergodic Markovian strategy. That is, the pair of strategies we mention there determines a Nash equilibrium of the Patroller in Uniform game. This is indeed a stronger result compared to what we were seeking for when initially considered a Stackelberg approach. Note that this Nash property holds in our later results as well (except for Section \ref{Section3.5}), though we make no claim here for the general existence of a Nash equilibrium in our game for all networks.
\end{remark}

\subsection{Attack duration $m$ odd}\label{Section3.2}
\indent

The solution to the uniformed Patroller game on the star is particularly
simple for $m$ odd, using an observation of the Associate Editor. Suppose $m=2j+1$ is odd. Suppose the Patroller adopts a random walk (always goes to a random adjacent node, never remains at the same node). In any $2j$ consecutive periods, he visits $j$ ends counting multiplicity, randomly chosen. The probability that a particular end $e$ is visited among these $j$ is given by $1-\left(\frac{n-1}{n}\right) ^{j}$ and so this is a lower bound on the interception probability. Similarly, if the Attacker begins her attack at any time the Patroller is away from her node, then the number $k$ of end nodes the Patroller visits during the rest of the attack satisfies $k\leq j,$ and since they are chosen randomly (independently) the probability that the Attacker's node is among them is $1-\left( \frac{n-1}{n}\right) ^{k}\leq-\left( \frac{n-1}{n}\right) ^{j}.$ Note that only the random walk always
gives $k=j.$ 

Hence, we have easily established the following:

\begin{proposition}\label{Proposition2}
Consider the Uniformed Patroller Problem on the star network $S_{n}$ with $n$ arcs and an attack difficulty of any odd integer $m.$ The unique optimal strategy for the Patroller is the random walk (which implies reflection at ends). An optimal strategy for the Attacker is an attack at any end node with an arbitrary delay. The value of the game is given by

\begin{equation}\label{eq8}
V=1-\left( \frac{n-1}{n}\right) ^{(m-1)/2}.
\end{equation}

\end{proposition}

\noindent We observe that this proves Conjecture \ref{Conjecture1} for all Star graphs $S_{n}$ when the attack duration $m$ is odd. 

\subsection{Attack duration $m=4$}\label{Section3.3}
\indent

Here, if the attack starts with the Patroller being at the center, then he gets two chances to intercept it, namely to find the correct end. However, if the attack starts with the Patroller being at an end, then he gets only one chance. In the first case, the Patroller can intercept the attack with any of the sequences 

\begin{equation*}
ce\_\_  \quad ,\quad cce\_ \quad ,\quad ccce \quad ,\quad ce^{'}\!\!ce,
\end{equation*}

\noindent where $e^{'}$ is any end node other than the attack end node $e$, while in the second case the Patroller can intercept the attack with either of the two sequences

\begin{equation*}
e^{'}\!\!ce\_ \quad ,\quad e^{'}\!\! cce\quad,\quad e^{'}\!\!e^{'}\!\!ce.
\end{equation*}

It follows that the attack will be intercepted with overall probability

\begin{gather}
\begin{aligned}\label{eq9}
a_{4}( n,p,s,q)&=q\cdot \bigl(1+r+r^{2}+( 1-p-r)\bigr)\cdot p+( 1-q)\cdot\bigl( s+s\cdot r+( 1-s)\cdot s\bigr)\cdot p
 \\&= -p\cdot ( r\cdot s-r^{2}-s^{2}+2\cdot s+p-2)\cdot q+p\cdot (
2\cdot s+r\cdot s-s^{2}),
\end{aligned}
\end{gather}

\noindent where $q$ is the conditional probability that the Patroller is at the center at the beginning of the attack given that he is not at the attack node $e$. Note that the coefficient of $q$ in \eqref{eq9} is given by the product of $p$ with the expression

\begin{gather*}
\begin{aligned}
( -r\cdot s+r^{2}+s^{2}-2\cdot s-p+2)&=2+( s^{2}-2\cdot s)-r\cdot s+r^{2}-p \\
\geq 2-1+-rs+r^{2}-p &\geq 1-\left( r+p\right) +r^{2}  \geq 1-\left( r+p\right)  \\
&\geq 0,
\end{aligned}
\end{gather*}

\noindent since $r+p\leq r+n\cdot p=1$, and $s\leq 1$. It follows that for fixed $n$ and $p$, the interception probability \eqref{eq9} is increasing in $q$. Thus, by the same reasoning that we have used for $m=2$ and $m=3$, it further follows that the Attacker should choose to wait for $d=2$ to attain $q=\hat{q}$, and the minimum interception probability of

\begin{equation}\label{eq10}
a_{4}=\frac{p}{1-p}\bigl(( p+r-1)\cdot s^{2}+(2-r-p\cdot r-r^{2}-2\cdot p)\cdot s+( 2\cdot r-p\cdot r+r^{3})\bigr) .
\end{equation}

To check that \eqref{eq10} is increasing in $s$, notice that the derivative with respect to $s$ in the bracketed quadratic above is given by

\begin{equation*}
2\cdot ( p+r-1)\cdot s+( 2-r-p\cdot r-r^{2}-2\cdot p) \geq 2\cdot ( p+r-1)\cdot 1+( 2-r-p\cdot r-r^{2}-2\cdot p) ,
\end{equation*}

\noindent since $r+p<1$, which is equivalent to 

\begin{equation*}
r\cdot (1- p-r) \geq 0,
\end{equation*}

\noindent since both factors are non-negative.

Consequently, it turns out that the Patroller maximizes the interception probability for fixed $p$ by choosing $s=1$ (reflecting at the ends). Taking $s=1$ in \eqref{eq10}, gives

\begin{gather}\label{eq11}
\begin{aligned}
a_{4}(n,p) &=\frac{p}{1-p}\cdot \bigl( 2\cdot r-p-2\cdot p\cdot r-r^{2}+r^{3}+1\bigr) 
 \\
&=\frac{-n^{3}\cdot p^{4}+2\cdot n^{2}\cdot p^{3}+2\cdot n\cdot
p^{3}-3\cdot n\cdot p^{2}-3\cdot p^{2}+3\cdot p}{1-p}.
\end{aligned}
\end{gather}

\begin{figure}[H]
  \begin{minipage}[b]{0.51\textwidth}
    \centering
   \includegraphics[scale=0.475]{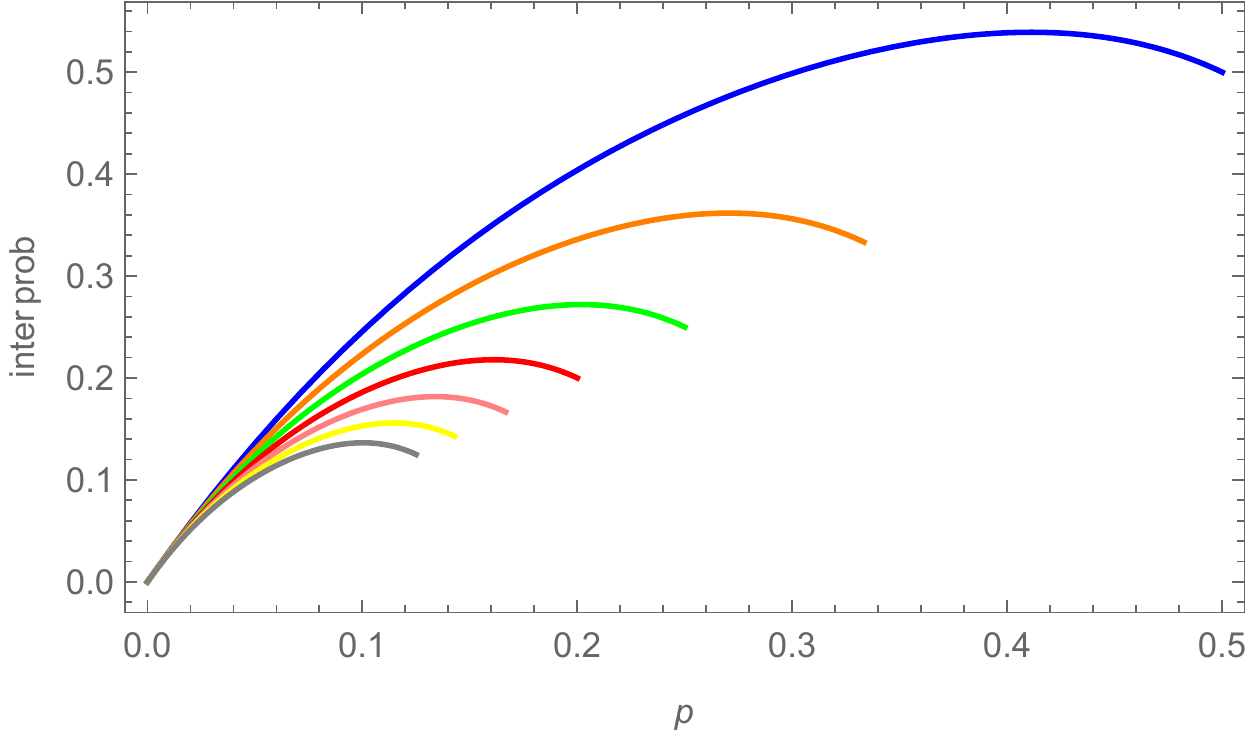}
    \captionof{figure}{$a_{4}(n,p)$ for $n=2$(blue), $\dots, 9$(black).}
\label{Figure7}
  \end{minipage}
\hfill
  \begin{minipage}[b]{0.43\textwidth}
    \centering
    \begin{tabular}{cccc}
	\hline\hline
      $n$ & $\hat{r}$ & $\hat{p}$ & $a_{4}(n,\hat{p})$ \\ \hline
        $2$ & $0.1778$ & $0.4111$ & $0.5391$ \\
        $3$ & $0.1885$ & $0.2705$ & $0.3618$ \\
        $4$ & $0.1924$ & $0.2019$ & $0.2720$ \\
        $5$ & $0.1945$ & $0.1611$ & $0.2179$ \\
        $6$ & $0.1960$ & $0.1340$ & $0.1817$ \\
        $7$ & $0.1971$ & $0.1147$ & $0.1559$ \\
        $8$ & $0.1976$ & $0.1003$ & $0.1364$ \\
        $9$ & $0.1981$ & $0.0891$ & $0.1213$ \\ \hline\hline
      \end{tabular}
      \captionof{table}{$a_{4}(n,\hat{p})$ for $n=2, \dots, 9$.}
\label{Table2}
    \end{minipage}
\end{figure}

In Figure \ref{Figure7} we plot the variation of the interception probability $a_{4}(n,p)$ with $p$, for different number of nodes, while in Table \ref{Table2} we give the maxima of $a_{4}(n,\hat{p})$ and the optimal $\hat{r}=1-n\cdot\hat{p}$.

\setlength{\tabcolsep}{10pt}
\renewcommand{\arraystretch}{0.95}

If we want, we can get an exact algebraic expression for the value $p=\hat{p}(n)$ which maximizes the interception probability \eqref{eq11}. We differentiate \eqref{eq11} with respect to $p$ and set the numerator of the resulting fraction $f(p)/(1-p)^{2}$ equal to $0$, where $f(p)$ is the fourth degree polynomial

\begin{equation*}
3-(6+6\cdot n)\cdot p +(3+9\cdot n+6\cot n^2)\cdot p^2+(-4\cdot n-4\cdot n^2+4\cdot n^3)\cdot p^3+(3\cdot n^3)^4.
\end{equation*}

The real root of this polynomial which is a probability is given by

\begin{equation*} 
\hat{p}(n)=\frac{2\cdot n^2+\frac{G}{\sqrt{2}}-3\cdot n^2\cdot\sqrt{\frac{C}{9\cdot N^4}+\frac{E}{3\cdot n^3\cdot D}+\frac{4\cdot \sqrt{2}\cdot F}{9\cdot n^2\cdot G}-\frac{D}{3\cdot n^3}}+2\cdot n+2}{6\cdot n^2},
\end{equation*}

\noindent where

\begin{align*}
\begin{gathered}
A=8\cdot n^6-18\cdot n^5+6\cdot n^4+9\cdot n^3-3\cdot n, \quad B=(n-1)^6\cdot n^3\cdot (32\cdot n^3+24\cdot n^2-3\cdot n-4),\\
D=\sqrt[3]{A+2\cdot\sqrt{B}-3\cdot n+1}, \quad C=8\cdot (n^2+n+1)^2-12\cdot n\cdot (2\cdot n^2+3\cdot n+1),\\
E=(4\cdot n^2-1)\cdot (n-1)^2, \quad G=\sqrt{C-\frac{6\cdot n\cdot E}{D}+6\cdot n\cdot D},\\
F=(4\cdot n^4+2\cdot n^3+6\cdot n^2+11\cdot n+4)\cdot(n-1)^2.
\end{gathered}
\end{align*}

\subsection{Asymptotic analysis of Star graph $S_{n}$ for $m=4$}\label{Section3.4}
\indent 

We now consider the optimal play on the star graph for difficulty $m=4$ when $n$ is large. Since $p$ goes to $0$, it is more transparent to work with the probability $r$ of staying at the center. We start by writing \eqref{eq11} in terms of $r$, using the fact that $p=(1-r)/n$, and calling it

\begin{equation*}
\pi(n,r)=\frac{(-1+r)\cdot \bigl(-1-r+2\cdot r^2+n\cdot(1+2\cdot r-r^2+r^3)\bigr)}{n\cdot(-1+n+r)},
\end{equation*}

\noindent so that,

\begin{equation*}
n\cdot \pi(n,r)\rightarrow poly(r)=r^4-2\cdot r^3+3\cdot r^2-r-1,\quad \text{as}\quad n\rightarrow \infty.
\end{equation*}

The first order condition on $r$ is given by the cubic

\begin{equation*}
4\cdot r^3-6\cdot r^2+6\cdot r-1=0,
\end{equation*}

\noindent with solution in $[0,1]$ of

\begin{equation*}
\begin{split}
\hat{r}_{\infty}&=\frac{1}{2}\Bigl(1-\bigl(-1+\sqrt{2}\bigr)^{-1/3}+\bigl(-1+\sqrt{2}\bigr)^{1/3}\Bigr)\simeq 0.20196.
\end{split}
\end{equation*}

We also have that

\begin{equation*}
\begin{split}
a\equiv poly(\hat{r}_{\infty})&=\frac{-3\Bigl(5-4\sqrt{2}-7\bigl(-1+\sqrt{2}\bigr)^{4/3}+\bigl(-1+\sqrt{2}\bigr)^{2/3}\cdot\bigl(-1+2\cdot \sqrt{2}\bigr)\Bigr)}{16\bigl(-1+\sqrt{2}\bigr)^{4/3}}\simeq 1.0944,
\end{split}
\end{equation*}

\noindent which implies that

\begin{equation*}
\pi(n,\hat{r}_{n})\rightarrow a/n \simeq 1.0944/n.
\end{equation*}

Hence, we show that for large $n$ the Patroller should stay at the center about $20\%$ of the time, which ensures him an interception probability of about $1.09/n$. The optimal delay remains~$2$.

\subsection{Giving the Patroller some memory}\label{Section3.5}
\indent

Our basic model gives the uniformed Patroller no memory, i.e. he is restricted to Markovian motion around the nodes of a graph. Here we modify the model for the star graph to allow the Patroller, when at the center, to remember if he has come there from another node (state $ec$) or has remained there from the previous period (state $cc$). The remaining state (denoted simply $e$) denotes his presence at an end other than the Attacker's node $E$. To simplify the analysis, we assume that when at an end (state $e$) he always reflects back to the center, as we found in Section \ref{Section3.1}. From the center he now has two probabilities, denoted $p$ and $s$, that he moves to a random end: the former $(p)$ from state $cc$ and the latter $(s)$ from state $ec$. If the Patroller follows such a strategy $(p,s)$, then conditional on not going to the Attacker's node $E$, his transition probabilities are given by the matrix

\begin{equation*}
\bordermatrix{
\mbox{states} & cc & ec & e \cr
cc & \frac{1-np}{1-p} & 0 & \frac{(n-1)p}{1-p} \cr
ec & \frac{1-ns}{1-s} & 0 & \frac{(n-1)s}{1-s} \cr
e & 0 & 1 & 0 \cr
}.
\end{equation*}

For simplicity, we will be mainly concerned with the case $n=3$ and $m=2$, where the corresponding matrix will be denoted by

\begin{equation}\label{eq12}
\renewcommand{\arraystretch}{1}
\setlength{\arraycolsep}{7.5pt}
A=\begin{pmatrix}
\frac{1-3p}{1-p} & 0 & \frac{2p}{1-p} \cr
\frac{1-3s}{1-s} & 0 & \frac{2s}{1-s} \cr
0 & 1 & 0 \cr
\end{pmatrix}.
\end{equation}

The Attacker knows that the Patroller's distribution over the three states, conditional on the fact that he has not been seen at the Attacker's chosen end node $E$, is given by the three tuple (with respective probabilities of being at states $cc$, $ec$ and $e$)

\begin{equation}\label{eq13}
x^{\kappa}=(0,1,0)\cdot A^{\kappa},
\end{equation}

\noindent after $\kappa$ periods away from $E$. We know that when he leaves $E$ he is at the center, and more specifically he is at the center having just arrived there form an end, so in state $ec$. That is, we have that $x^{1}=(0,1,0)$. If the attack takes place with delay $d$, and $x^{d}=(x^{d}_{1},x^{d}_{2},x^{d}_{3})$, it will be intercepted if and only if the Patroller moves to $E$ in the next period, which has probability
\begin{equation}\label{eq14}
P(p,s,d)=p\cdot x^{d}_{1}+s\cdot x^{d}_{2}=x^{d}\cdot (p,s,0).
\end{equation}

Thus, in our Stackelberg approach, the Patroller's objective is to choose $p$ and $s$ (both in the interval $[0,1/3]$) to maximize the minimum of \eqref{eq14} over $d$ (chosen by the Attacker). So we have that 
\begin{equation}\label{eq15}
V=\max_{p,s\in[0,1/3]}\min_{1\leq d\leq D} x^{d}\cdot (p,s,0).
\end{equation}

\noindent Note that $x^{d}$ is also a function of $p$ and $s$, though this is not indicated in our notation.

We solve \eqref{eq15} by first numerical and then exact algebraic methods. For our numerical work we take $D=10$. Unlike the similar problem in Section \ref{Section3.1}, it turns out that there is not a unique delay $d$ which the Attacker can adopt without knowing $p$ and $s$ (here it is either $d=2$ or $d=3$). If the Attacker moves first, then she will have to use a mixture of $d=2$ and $3$, but even then she cannot reduce the interception probability $P$ to the maximum value $V$. Numerical work fixes $p$ and $s$, generates the sequence of interception probabilities $x^{d}\cdot (p,s,0)$ in \eqref{eq15}, and selects the delay $d$ which minimizes the entry. By varying $p$ and $s$ in a grid of values, we find that the optimal $p$ is near $0.3$, the optimal $s$ is near $0.2$, and the minimizing delay $d$ is either $2$ or $3$. To illustrate these ideas we give the sequences $x^{d}\cdot (p,s,0)$ for the two strategy pairs $p=0.30$, $s=0.21$ and $p=0.30$, $s=0.22$ as

\begin{equation}
\begin{split}
0.21,\,0.141,\,0.132^{*},\,0.162,\,0.139,\,0.149,\,0.148,\,0.145,\,0.148,\,0.147,\quad &\text{for} \, (p,s)=(0.30,0.21),
\\
0.22,\,0.131^{*},\,0.143,\,0.159,\,0.142,\,0.152,\,0.148,\,0.148,\,0.149,\,0.148,\quad &\text{for} \, (p,s)=(0.30,0.22).
\end{split}
\end{equation}

\noindent If the Patroller adopts $(0.30,0.21)$, then the optimal Attacker's response is $d=3$ with interception probability $0.132$, while if the Patroller adopts $(0.30,0.22)$, then the optimal response is $d=2$ with interception probability $0.131$. Thus we have a numerical function $h(p,s)$, which is calculated as $h(0.30,0.21)=0.132$ and $h(0.20,0.22)=0.131$. By carrying out such numerical work over $p$ and $s$ and taking the maximum, a good approximation to the maximizing values can be found (similar to numerical methods we implement in Sections \ref{Section4}, \ref{Section5} and \ref{Section6}).

To get more accurate results we adopt algebraic methods, based on our numerical results for $d$. As functions of $p$, $s$, the first three terms of the $x^{\kappa}$ distribution sequence \eqref{eq13} are given respectively by

\begin{equation*}
\begin{split}
(0,1,0),\qquad (0,1,0)\cdot A&=\biggl(\frac{1-3s}{1-s},0,\frac{2s}{1-s}\biggr),\qquad \text{and}
\\
\biggl(\frac{1-3s}{1-s},0,\frac{2s}{1-s}\biggr)A&=\biggl(\frac{(1-3p)(1-3s)}{(1-p)(1-s)},\frac{4s}{1-s},\frac{2p(1-3s)}{(1-p)(1-s)}\biggr).
\end{split}
\end{equation*}

\noindent Then the first terms of the sequence $x^{\kappa}\cdot (p,s,0)$ of interception probabilities are

\begin{equation*}
\begin{split}
(0,1,0)\cdot (p,s,0)&=s,\\
\biggl(\frac{1-3s}{1-s},0,\frac{2s}{1-s}\biggr)\cdot (p,s,0)&=\frac{p(1-3s)}{1-s}\equiv u_{2}(p,s),\\
\biggl(\frac{(1-3p)(1-3s)}{(1-p)(1-s)},\frac{4s}{1-s},\frac{2p(1-3s)}{(1-p)(1-s)}\biggr)\cdot A&=\frac{p(1-3p)(1-3s)}{(1-p)(1-s)}+\frac{2s^{2}}{1-s}\equiv u_{3}(p,s),
\end{split}
\end{equation*}

\noindent where $u_{2}(p,s)$ and $u_{3}(p,s)$ are the payoffs if the Attacker chooses $d=2$ or $d=3$ while the Patroller chooses his two probabilities of moving from the center as $p$ and $s$, as described above. Thus, given our knowledge from numerical work that the optimal waiting time $d$ is always $2$ or $3$, we can rewrite the value problem \eqref{eq15} as 

\begin{equation*}
V=\max_{p,s\in [0,1/3]}\min (u_{2}(p,s),u_{3}(p,s)).
\end{equation*}

Solving the equation $u_{2}(p,s)=u_{3}(p,s)$, we get the optimal $s$ as a function of $p$

\begin{equation*}
\hat{s}(p)=\frac{-3p^{2}+\sqrt{4p^2-4p^3+9p^4}}{2(1-p)}.
\end{equation*}

\noindent It follows that the value of the game is given by $\max_{p}u(p)$, where $u(p)=u_{2}(p,\hat{s}(p))=u_{3}(p,\hat{s}(p))$ is given by the common formula

\begin{equation}\label{eq16}
u(p)=\frac{p\bigl(2-2p+9p^2-3p\sqrt{4-4p+9p^2}\bigr)}{2-2p+3p^2-p\sqrt{4-4p+9p^2}}.
\end{equation}

\noindent The maximum of function \eqref{eq16} on $[0,1]$ cannot be simply expressed, but it is easy to approximate it at $\bar{p}\simeq 0.305$, with $\bar{s}=\hat{s}(\bar{p})\simeq 0.217$ and $V=u(\bar{p})\simeq 0.136$. Recall that the similar value (interception probability) when the Patroller has no memory was shown in Section \ref{Section3.1} to be $5-2\sqrt{6}\simeq 0.101$. Thus a Patroller with limited memory intercepts the attack about $13\%$ of the time compared with $10\%$ for a Patroller with no memory, which is about a $34\%$ increase. It is also useful to compare the optimal probabilities of staying still when at the center. For the Markovian Patroller of Section \ref{Section3.1} this probability (i.e. the $\hat{r}$ of \eqref{eq6}) is $0.449$. For the Patroller with memory, it is $1-3\bar{p}\simeq 0.0855$ if he was previously again at the center, and $1-3\bar{s}\simeq 0.349$ if he has just come there from an end.

Hence, we have the following result:

\begin{proposition}\label{Proposition3}
Consider that the Patroller, when at the center, knows whether he has just arrived or was already there the previous period. Let $p$ and $s$ denote his choice of probability to move to an end from the center, respectively, in the two cases. For the Stackelberg problem, the Patroller maximizes the interception probability at $0.136$ by choosing $p$ and $s$ with respective probabilities approximately $0.305$ and $0.136$. The Patroller's interception probability is an increase of approximately $36\%$ over what he had without any memory. At the optimal values of $p$ and $q$, both $d=2$ and $d=3$ are optimal for the Attacker.
\end{proposition}

\section{The Line Network $L_{n}$}\label{Section4}

The line network $L_{n}$ consists of $n$ nodes connected by $n-1$ edges in a line. We label the nodes $1,2,\dots,n$ from one end to the other. We restrict the Patroller to ergodic Markovian strategies, introducing the following setting. At the ends he reflects with probability $\kappa$, while he stays put with probability $\ell=1-\kappa$. From any other internal node $j\neq \mathit{1,n}$, he moves towards the closest end with probability $p_{j}$, he moves towards the furthest end with probability $q_{j}$, while he remains at the node with probability $r_{j}=1-p_{j}-q_{j}$. To be consistent with this notation, we set $p_{0}=p_{n}=0$, $q_{0}=q_{n}=\kappa$, and $r_{0}=r_{n}=\ell$. In the special case when $L_{n}$ consists of an odd number of nodes, the Patroller shifts from the center towards the two ends with probability $c$, staying still with probability $s=1-2\cdot c$.

We further assume that the Patroller is symmetric in his random patrolling, namely that 

\begin{equation*}
p_{j}=p_{n+1-j}\quad and\quad q_{j}=q_{n+1-j}.
\end{equation*}

\noindent Then, the transition matrix $A_{n}$ characterizing the Patroller's walk consists of $n-1$ parameters. 

This last assumption establishes a symmetry for the Attacker's strategy as well. Provided that the Patroller `announces' an ergodic Markovian  patrol $(p_{i},q_{i},c,\kappa)_{i\in K}$, the Attacker expects equivalent interception probabilities (i.e. equivalent payoffs) for her attack adopting either of the strategies $(i,d)$, $(n+1-i,d)$, where $i$, $n+1-i$ are symmetric, with respect to the center, nodes, and $d$ is the delay.

We focus our attention on the analysis of an attack planned to take place at an end node, in particular at node $\mathit{1}$. The Patroller's initial position distribution in this case is $x^{(0)}=(1,0,\dots,0)$. His future position distribution conditional to him not having returned to the attack node $\mathit{1}$ for $t$ consecutive periods is denoted by $x^{(t)}=(x_{1}^{(t)},\dots,x_{n}^{(t)})$, where $x_{1}^{(t)}=0$ for $t\geq 1$. Notice that whenever the Patroller returns to the attack node, $t$ is reset to zero and his position distribution becomes again $x^{(0)}$. For that, we refer to $x^{(t)}$ as \textit{the Patroller's away distribution} at $t\geq 1$. For example, the first two periods that the Patroller is away from node $\mathit{1}$, his away distribution is

\begin{equation*}
x^{(1)}=(0,1,0,\dots,0),\;\; x^{(2)}=(0,\frac{r_{2}}{1-p_{2}},\frac{q_{2}}{1-p_{2}},0,\dots,0),\;\; etc.
\end{equation*}

More generally, if we define the row vector $y^{(t)}$, for $t\geq 0$, by the equation

\begin{equation*}
y^{(t)}=x^{(t)}\times A_{n},
\end{equation*}

\noindent then the Patroller's away distribution at $t\geq 1$ is given by the following iteration formula

\begin{equation}\label{eq18}
x^{(t)}=\frac{x^{(t-1)}\times A_{n}-(y_{1}^{(t-1)},0,\dots,0)}{1-y_{1}^{(t-1)}}.
\end{equation}

Since we adopt a Stackelberg approach where the Patroller is the leader, in fact we let the Attacker observe the Patroller's walk on the network prior to deciding her strategy. The Patroller's objective is to `announce' the optimum strategy $(\hat{p}_{i},\hat{q}_{i},\hat{\kappa},\hat{c})$ that maximizes the interception probability of an attack of known difficulty $m$ for what is rationally the Attacker's delay $\hat{d}$, while the Attacker's objective is to decide the optimum delay $\hat{d}$ that minimizes the interception probability of her attack for what is the announced-observed Markovian patrol. Hence,  the value $V$ of the game, is given by

\begin{equation}\label{eq19}
V=\max_{A_{n}}\min_{d}\pi_{m}(p_{i},q_{i},\kappa,c,d)=\pi_{m}(\hat{p}_{i},\hat{q}_{i},\hat{\kappa},\hat{c},\hat{d}),
\end{equation}

\noindent where $\pi_{m}(p_{i},q_{i},\kappa,c,d)$ is the interception probability of an attack at node $\textit{1}$ of duration $m$, under the patrol $A_{n}$ and with delay $d$. We solve the game numerically for the line networks $L_{4}$, $L_{5}$, for given values of $m$. $L_{3}$ has been examined as a sub-case of the star network; it is equivalent to $S_{2}$.

\subsection{The Network $L_{4}$}\label{Section4.1}
\indent

First we analyze the uniformed patrolling game on the line network $L_{4}$ drawn below.
\begin{figure}[h!]
\centering
 \includegraphics[scale=0.85]{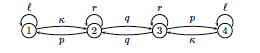}
\caption{The Line network $L_{4}$.}
\label{Figure8}
\end{figure}

To avoid unnecessary subscripts we introduce notation $p_{2}=p$, $q_{2}=q$, $r_{2}=r$. The general vector formula \eqref{eq18} for recursively calculating Patroller's away distribution, reduces to the system

\begin{gather}
\begin{aligned}\label{eq20}
(1-p\cdot x_{2}^{(t-1)})\cdot x_{2}^{(t)} &=(1-p-q)\cdot x_{2}^{(t-1)}+q\cdot x_{3}^{(t-1)},\\
(1-p\cdot x_{2}^{(t-1)})\cdot x_{3}^{(t)} &=(q-\kappa)\cdot x_{2}^{(t-1)}+(1-p-q-\kappa)\cdot x_{3}^{(t-1)}+\kappa.
\end{aligned}
\end{gather}

For fixed $p$, $q$ and $\kappa$, system \eqref{eq20} defines a continuously differentiable mapping $f_{2}:\Delta^{2}\mapsto\Delta^{2}$. Any vector of the form $(0,x_{2},x_{3},x_{4})$, such that $x_{2},x_{3},x_{4}\in (0,1)$ and $x_{2}+x_{3}+x_{4}=1$, will be a fixed point of $f_{2}$, if it also satisfies the system of equations

\begin{gather}
\begin{aligned}\label{eq21}
(p+q-p\cdot x_{2})\cdot x_{2} &=q\cdot x_{3},\\
(\kappa+p+q-p\cdot x_{2})\cdot x_{3} &=\kappa+(q-\kappa)\cdot x_{2}.
\end{aligned}
\end{gather}

\noindent For the optimal patrolling strategies we estimate below, system \eqref{eq21} has a unique solution, which gives the limiting Patroller's away distribution under these optimal patrols (see Table \ref{Table3}).  

We calculate the interception probability $\pi_{m}$ of an attack at node $\mathit{1}$, considering five cases $m=2,3,4,5,6$. If the Attacker waits for $t=d$ time periods before she attacks, then we can compute the interception probability by conditioning on the Patroller's location at time $t=d$.

\begin{align*}
\pi_{2}(p,q,\kappa,d) &=p\cdot x_{2}^{(d)},\\
\pi_{3}(p,q,\kappa,d) &=p\cdot\bigl(1+(1-p-q)\bigr)\cdot x_{2}^{(d)}+p\cdot q\cdot x_{3}^{(d)},\\
\begin{split}
\pi_{4}(p,q,\kappa,d) &=p\cdot\bigl(1+(1-p-q)+(1-p-q)^2+q^2\bigr)\cdot x_{2}^{(d)} +p\cdot q\cdot \bigl(1+2\cdot (1-p-q)\bigr) \cdot x_{3}^{(d)}\\&+p\cdot q\cdot\kappa\cdot x_{4}^{(d)},
\end{split}
\end{align*}

\noindent whereas we omit the lengthy explicit formulas for $\pi_{5}(p,q,\kappa,d)$ and $\pi_{6}(p,q,\kappa,d)$.

We estimate numerically, for $D=15$, the optimal delay, the optimal patrol, and the corresponding interception probability for the above five cases, and we present our findings in Table \ref{Table3} rounded to four decimal places. Note that we include an extra column with the interception probabilities under the Patroller's optimal walk $(\hat{p},\hat{q},\hat{\kappa})$, but for delay $d\gg 15$. One can think of $\pi_{m}(\hat{p},\hat{q},\hat{\kappa},\infty)$ as the interception probability the Attacker should expect if she waits for a long time before attacking.

\setlength{\tabcolsep}{13pt}
\renewcommand{\arraystretch}{0.8}
\begin{table}[h!]
\centering
\begin{tabular}{ ccccc } 
\hline\hline
 for $D=15$    & $(\hat{p},\hat{q},\hat{\kappa})$ & $\hat{d}$ & $\pi_{m}(\hat{p},\hat{q},\hat{\kappa},\hat{d})$ & $\pi_{m}(\hat{p},\hat{q},\hat{\kappa},\infty)$   \\ \hline
 $m=2$ & $(0.3935,0.3309,1)$ & $4$ & 0.1032 & $0.1067$ \\ 
 $m=3$ & $(0.5,0.5,1)$ & $2$ & 0.25 & $0.25$  \\ 
 $m=4$ & $(0.4317,0.4076,1)$ & $4$ & $0.2960$ & $0.3158$ \\
 $m=5$ & $(0.5,0.5,1)$ & $2$ & $0.4375$ & $0.4375$ \\ 
 $m=6$ & $(0.4974,0.4267,1)$ & $4$ & $0.4551$ & $0.4766$ \\ \hline\hline
\end{tabular}
\caption{Optimal game values for $L_{4}$ (Attack node $\mathit{1}$)}
\label{Table3}
\end{table}

\begin{remark}\label{Remark2}
We find that the Patroller should always reflect at the ends ($\hat{\kappa}=1$). Additionally, for an odd attack duration ($m=3,5$) the optimal patrol is a random walk ($\hat{p}=\hat{q}=0.5$) that reflects at the boundaries, that is, the Patroller should never remain at the same node for two consecutive periods. Accordingly, the Attacker's best strategy is $(\mathit{1},2)$. For an even attack duration ($m=2,4,6$) the Patroller optimally remains at an internal node with non-zero probability, which, however, decreases as we increase $m$ ($r_{m=2}>r_{m=4}>r_{m=6}$). We further observe that $\hat{p}>\hat{q}$ for an even $m$. For example, for $m=4$ the optimal $p$ and $q$ are $0.4317$ and $0.4076$ respectively. The Attacker's optimal response to this Markovian patrol is to initiate her attack, at node $\mathit{1}$, in the fourth period ($d=4$) that the Patroller is away. Under these strategies the attack is intercepted with probability $0.2960$. In the alternative case when the Attacker chooses to wait for a very long time at the Patroller's absence before attacking ($d\gg 1$), her attack will be more likely to be intercepted (with probability $0.3158$). 
\end{remark}

However, there are two issues that naturally arise regarding our analysis, and need to be further investigated. To treat both of them we adopt a qualitative approach.

$(i)$\textit{How confident we can be that for a maximum waiting time $D>15$ the Attacker cannot achieve a lower interception probability?} 

\begin{figure}[h!]
\centering
\subfigure[$m=2$]{\includegraphics[scale=0.41]{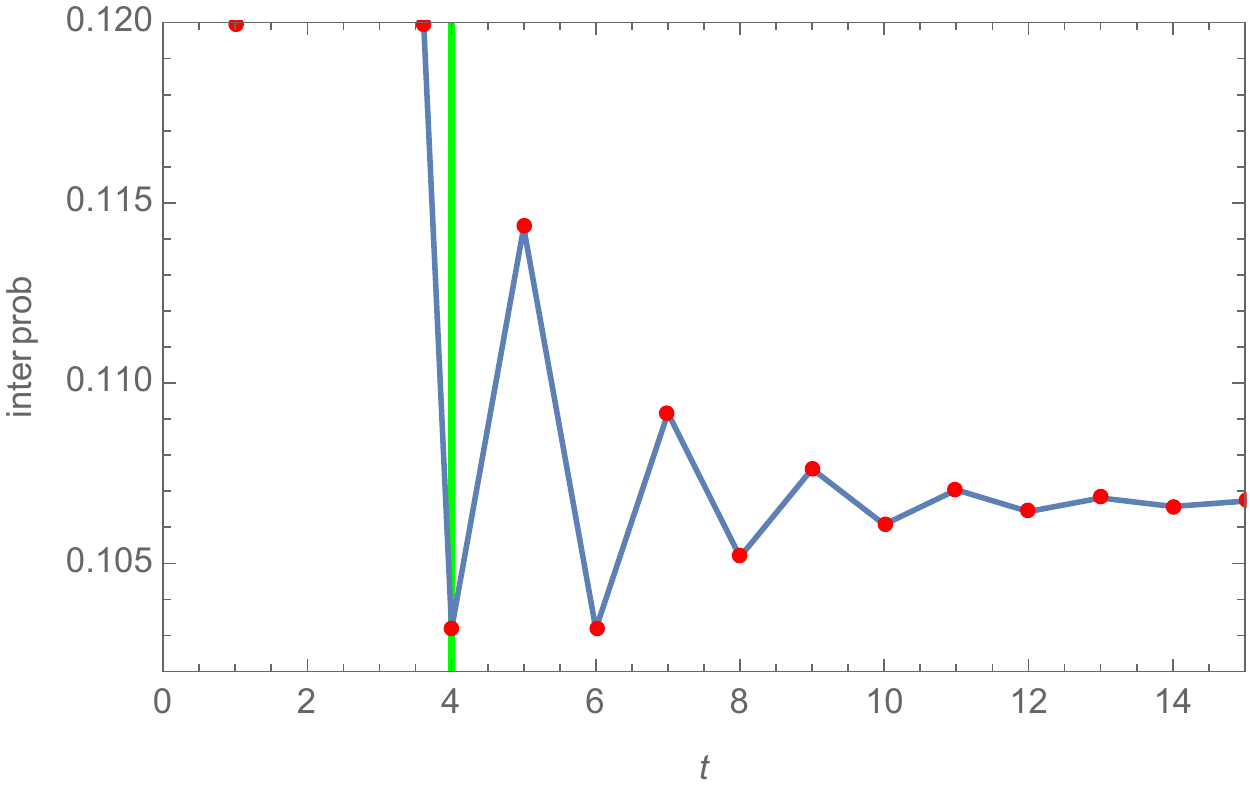}}
\subfigure[$m=4$]{\includegraphics[scale=0.41]{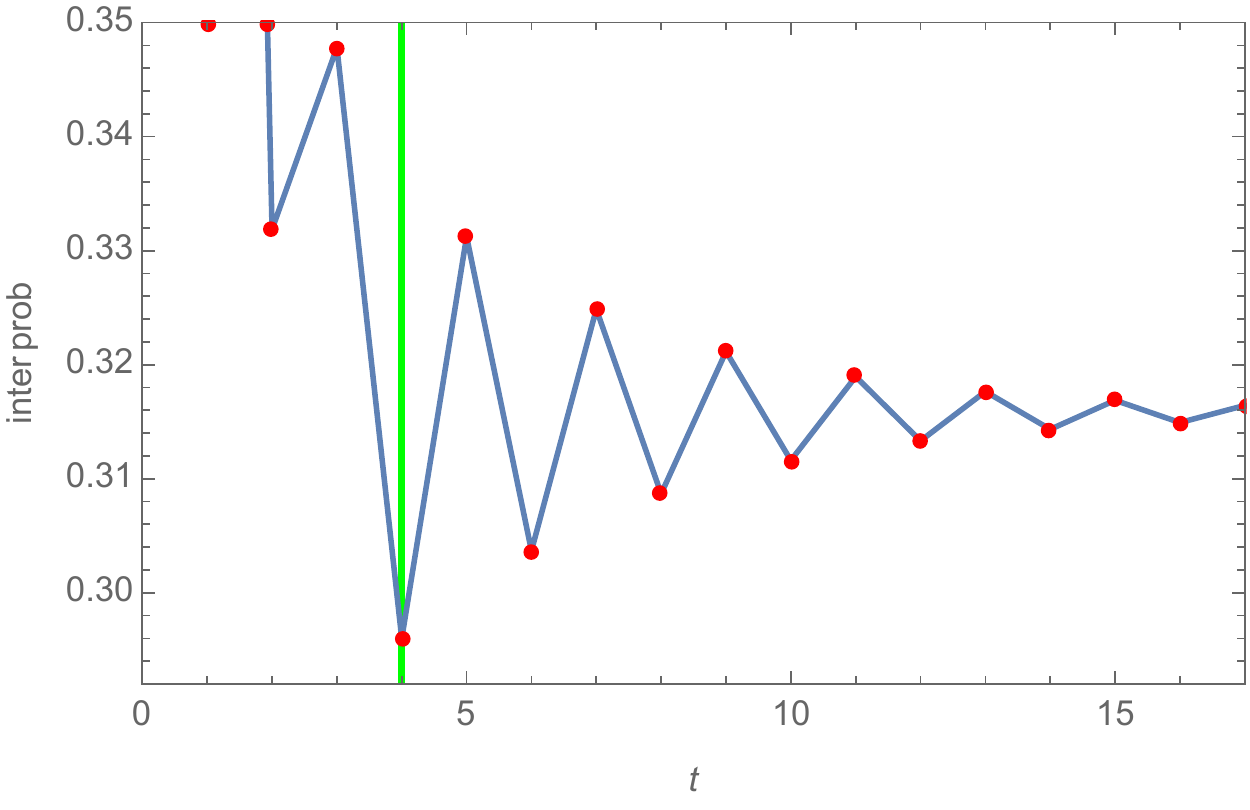}}
\subfigure[$m=6$]{\includegraphics[scale=0.41]{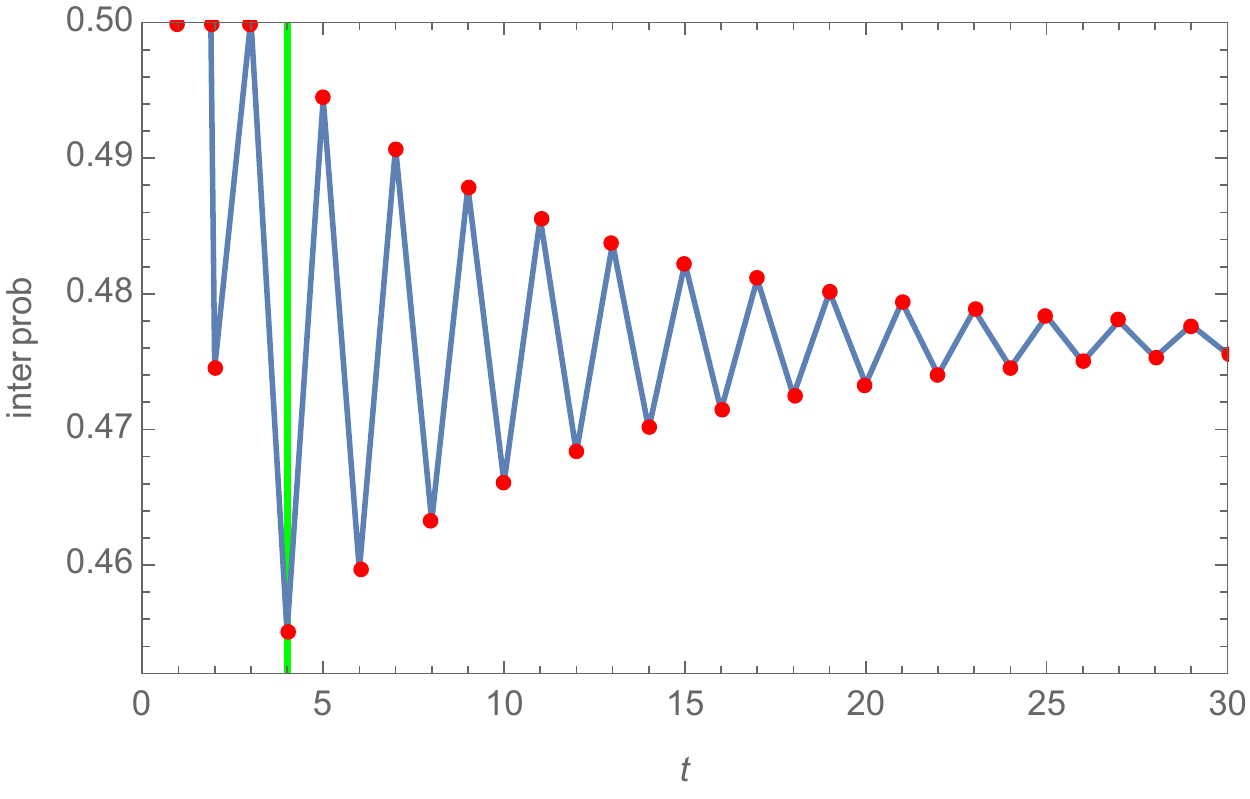}}
\caption{Interception probability for an attack at node $\mathit{1}$ of $L_{4}$, under the optimal patrol, for increasing delay.}
\label{Figure9}
\end{figure}

\noindent To investigate the first issue, we generate the interception probabilities $\pi_{m}(p,q,\kappa,d)$ for increasing values of $d$ and for the optimal $\hat{p}, \hat{q}, \hat{\kappa}$ given in Table \ref{Table3}, and check if we reach the limiting interception probabilities $\pi_{m}(\hat{p},\hat{q},\hat{\kappa},\infty)$ without crossing below the optimal values $\pi_{m}(\hat{p},\hat{q},\hat{\kappa},\hat{d})$ that we get for $D=15$ (the vertical green line indicates the optimal delay estimated for $D=15$). As seen in Figure \ref{Figure9}, this appears to be the case for $m=2$(\ref{Figure9}a),$4$(\ref{Figure9}b),$6$(\ref{Figure9}c), while for $m=3,5$ we already reach the Patroller's stationary away distribution from the second~period.

$(ii)$\textit{Can the Attacker do any better by attacking a penultimate node instead of an end node?}

\noindent Regarding the second issue, for each value of $m$ we consider the optimal patrol $(\hat{p},\hat{q},\hat{\kappa})$ that we estimate for an attack at node $\mathit{1}$ (e.g. $(0.4947,0.4267,1)$ for $m=6$), and generate the analogous interception probability for an attack at node $\mathit{2}$ for increasing values of $d$. For $m=6$ we get 

\begin{figure}[H]
\centering
\includegraphics[scale=0.41]{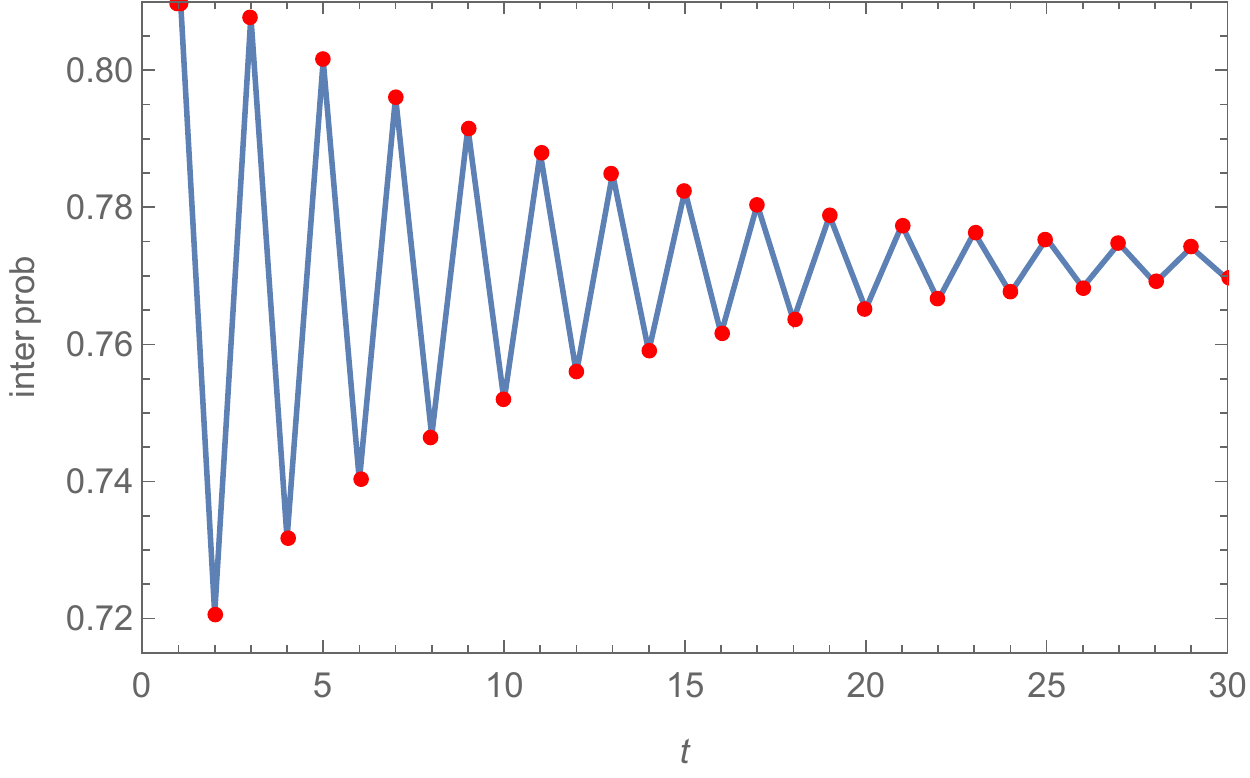}
\caption{Interception probability for an attack at node $\mathit{2}$ of $L_{4}$, under the optimal patrol estimated for an attack at node $\mathit{1}$, for increasing delay, for $m=6$.}
\label{Figure10}
\end{figure}

\noindent As seen in Figure \ref{Figure10}, the minimum interception probability ($0.7207$) for the attack strategy $(\mathit{2},1)$ is greater than the corresponding probability ($0.4551$) for the attack strategy $(\mathit{1},4)$ (both estimated under the optimal patrol $(0.4974,0.4267,1)$, see Table \ref{Table3}). Hence, attacking at an end with $d=4$ is optimal. Similar results we obtain for the rest of the cases, as given in Table~\ref{Table4}.

\setlength{\tabcolsep}{13pt}
\renewcommand{\arraystretch}{0.85}
\begin{table}[H]
\centering
\begin{tabular}{ cccc } 
\hline\hline
 for $(\hat{p},\hat{q},\hat{\kappa})$  & $\pi_{m}$ at node $\mathit{1}$, $(i,\hat{d})$ & $\pi_{m}$ at node $\mathit{2}$, $(i,\hat{d})$   \\ \hline
 $m=2$ & $0.1032$, $(\mathit{1},4)$ & $(\mathit{2},2)-0.1363$ \\ 
 $m=3$ & $0.2500$, $(\mathit{1},2)$ & $(\mathit{2},1)-0.5000$  \\ 
 $m=4$ & $0.2960$, $(\mathit{1},4)$ & $(\mathit{2},2)-0.5237$ \\
 $m=5$ & $0.4375$, $(\mathit{1},2)$ & $(\mathit{2},1)-0.7500$ \\ 
 $m=6$ & $0.4551$, $(\mathit{1},4)$ & $(\mathit{2},2)-0.7207$ \\ \hline\hline
\end{tabular}
\caption{Compare attacks at nodes $\mathit{1}$ and $\mathit{2}$ under the optimal patrol for node $\mathit{1}$}
\label{Table4}
\end{table}

\subsection{The Network $L_{5}$}\label{Section4.2}
\indent 

Next we consider the line network $L_{5}$ which differs from $L_{4}$ in having a central node, node~$\mathit{3}$.

\begin{figure}[H]
\centering
\includegraphics[scale=0.85]{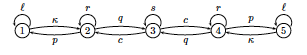}
\caption{The Line network $L_{5}$.}
\label{Figure11}
\end{figure}

Expanding the general recursive formula \eqref{eq18}, we can rewrite the Patroller's away distribution for $L_{5}$ into the following system of equations 

\begin{gather}
\begin{aligned}\label{eq23}
(1-p\cdot x_{2}^{(t-1)})\cdot x_{2}^{(t)} &=(1-p-q)\cdot x_{2}^{(t-1)}+c\cdot x_{3}^{(t-1)},\\
(1-p\cdot x_{2}^{(t-1)})\cdot x_{3}^{(t)} &=q\cdot x_{2}^{(t-1)}+(1-2\cdot c)\cdot x_{3}^{(t-1)}+q\cdot x_{4}^{(t-1)},\\
(1-p\cdot x_{2}^{(t-1)})\cdot x_{4}^{(t)} &=\kappa-\kappa\cdot x_{2}^{(t-1)}+(c-\kappa)\cdot x_{3}^{(t-1)}+(1-p-q-\kappa)\cdot x_{4}^{(t-1)}.
\end{aligned}
\end{gather}

For fixed $p$, $q$, $c$, $\kappa$, system \eqref{eq23} defines a continuously differentiable mapping $f_{3}:\Delta^{3}\mapsto\Delta^{3}$. To determine the limiting behavior of the Patroller's away distribution, we note that any vector of the form $(0,x_{2},x_{3},x_{4},x_{5})$, such that $x_{2},x_{3},x_{4},x_{5}\in (0,1)$ and $x_{2}+x_{3}+x_{4}+x_{5}=1$, is a fixed point of $f_{3}$, if it also satisfies the system 

\begin{gather}
\begin{aligned}\label{eq24}
(p+q-p\cdot x_{2})\cdot x_{2} &=c\cdot x_{3},\\
(2\cdot c-p\cdot x_{2})\cdot x_{3} &=q\cdot (x_{2}+ x_{4}),\\
(p+q+\kappa-p\cdot x_{2})\cdot x_{4} &=\kappa-\kappa\cdot x_{2}+(c-\kappa)\cdot x_{3}.
\end{aligned}
\end{gather}

\noindent For the optimal patrolling strategies we estimate below, system \eqref{eq24} has a unique solution, that is the Patroller's away distribution if the Attacker waits long at his absence before attacking.

We determine the interception probability $\pi_{m}(p,q,c,\kappa,d)$ for the Patroller's Markovian parameters $p,q,c,\kappa$, and the Attacker's delay $d$.  We consider the same five cases, $m=2,3,4,5,6$, and respectively we have

\begin{align*}
\pi_{2}(p,q,c,\kappa,d)&=p\cdot x_{2}^{(d)},\\
\pi_{3}(p,q,c,\kappa,d)&=p\cdot \bigl(1+(1-q-p)\bigr)\cdot x^{(d)}_{2}+p\cdot c\cdot x^{(d)}_{3},\\
\begin{split}
\pi_{4}(p,q,c,\kappa,d)&=p\cdot \bigl(1+(1-p-q)+q\cdot c+(1-p-q)^{2}\bigr)\cdot x_{2}^{(d)}\\&+c\cdot p\cdot \bigl(1+(1-p-q)+(1-2\cdot c)\bigr)\cdot x_{3}^{(d)}+ q\cdot c \cdot p\cdot x_{4}^{(d)},
\end{split}
\end{align*}

\noindent where we omit the explicit formulas for $\pi_{5}$ and $\pi_{6}$. Likewise, we estimate numerically, for $D=15$, the optimal game values presented in Table \ref{Table5}, rounded to four decimal places.

\setlength{\tabcolsep}{13pt}
\renewcommand{\arraystretch}{0.8}
\begin{table}[h!]
\centering
\begin{tabular}{ ccccc } 
\hline\hline
 for $D=15$    & $(\hat{p},\hat{q},\hat{c},\hat{\kappa})$ & $\hat{d}$ & $\pi_{m}(\hat{p},\hat{q},\hat{c},\hat{\kappa},\hat{d})$ & $\pi_{m}(\hat{p},\hat{q},\hat{c},\hat{\kappa},\infty)$   \\ \hline
 $m=2$ & $(0.4342,0.4110,0.4096,1)$ & $10$ & $0.0646$ & $0.0664$  \\ 
 $m=3$ & $(0.5,0.5,0.5,1)$ & $12$ & $0.1464$ & $0.1464$  \\ 
 $m=4$ & $(0.4663,0.4330,0.4216,1)$ & $8$ & $0.1890$ & $0.1932$  \\ 
 $m=5$ & $(0.5,0.5,0.5,1)$ & $14$ & $0.2714$ & $0.2714$  \\ 
 $m=6$ & $(0.4821,0.4598,0.4332,1)$ & $8$ & $0.3$ & $0.3088$\\ \hline\hline
\end{tabular}
\caption{Optimal game values for $L_{5}$ (Attack node $\mathit{1}$)}
\label{Table5}
\end{table}

\begin{remark}\label{Remark3}
As expected from the previous analysis for $L_{4}$, the Patroller should  reflect at the ends in this case too. Similarly, from Table \ref{Table5} we conclude that a random walk is the optimal patrol when the duration of the attack is odd (i.e. for $m=3,5$). A third observation is that for an even attack duration (i.e. for $m=2,4,6$) the optimal $\hat{p}$ is always greater than the optimal $\hat{q}$, namely the optimal strategy for the Patroller is to move towards the closest end with higher probability than move towards the center. We see as well that increasing the attack duration for even values, the optimal patrol converges to a random walk, since $\hat{p}$, $\hat{q}$ and $\hat{c}$ increase in even~$m$. Our findings are in principle identical with those obtained in subsection \ref{Section4.1} for the network $L_{4}$.
\end{remark}

\begin{figure}[h!]
\centering
\subfigure[$m=2$]{\includegraphics[scale=0.41]{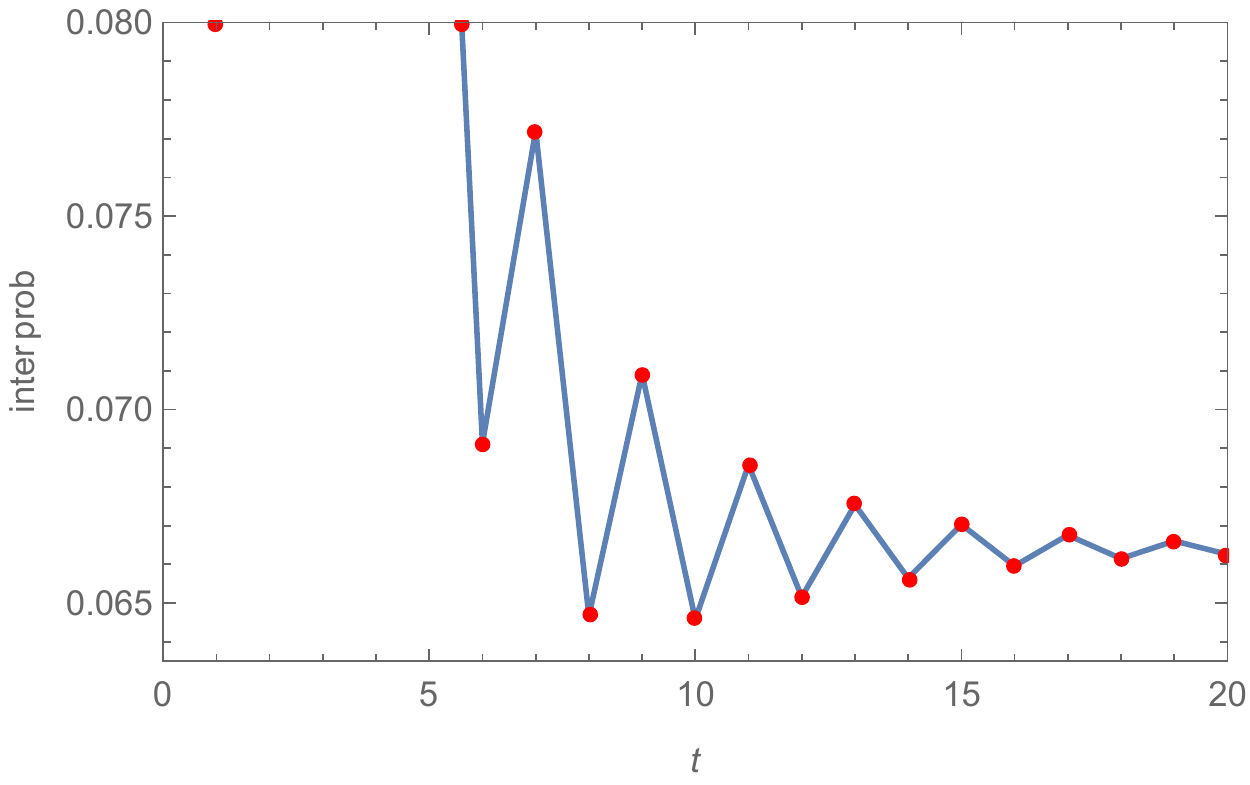}}
\subfigure[$m=4$]{\includegraphics[scale=0.41]{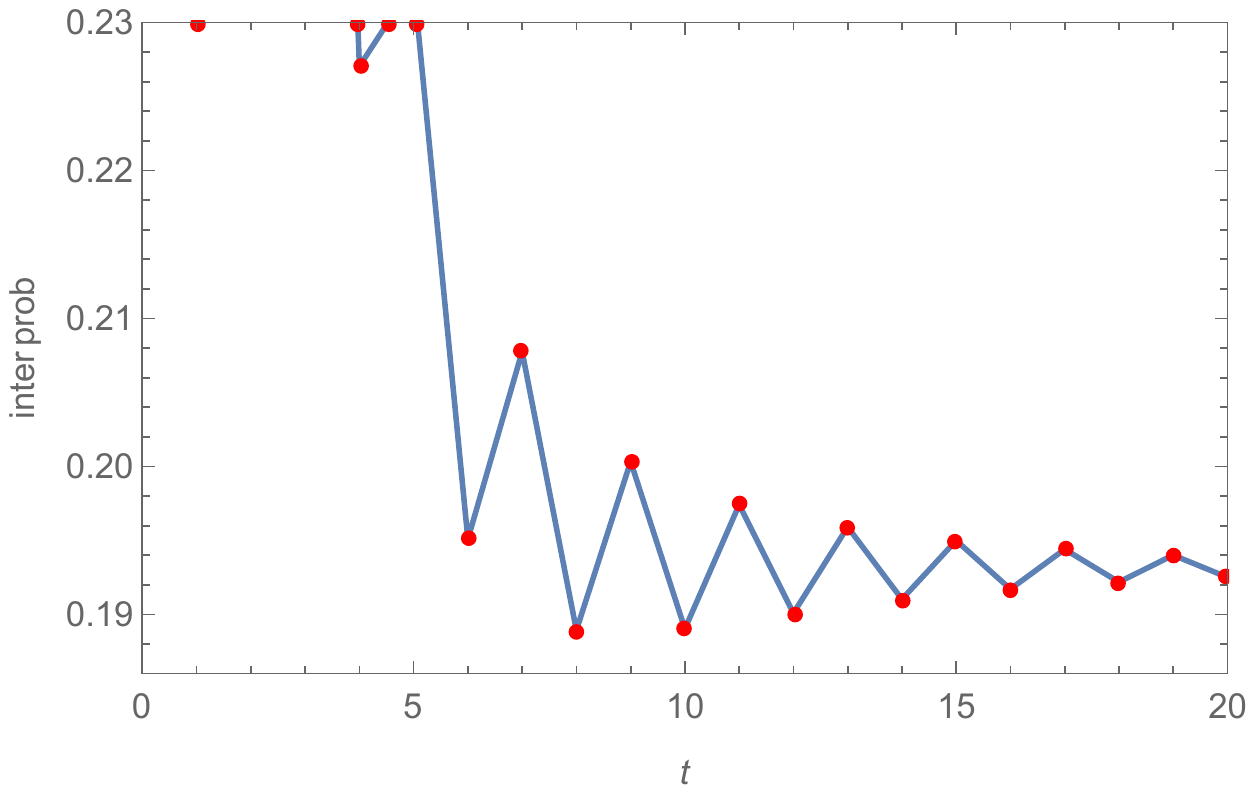}}
\subfigure[$m=6$]{\includegraphics[scale=0.41]{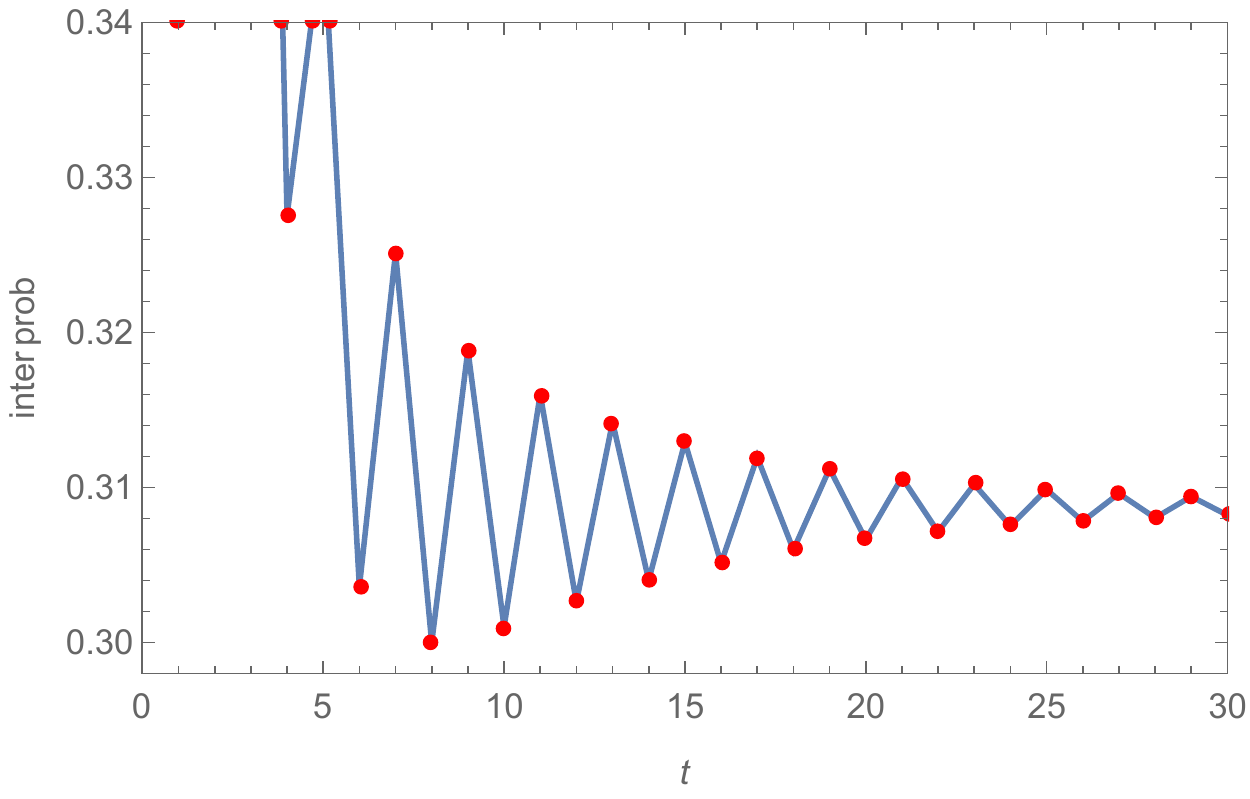}}
\subfigure[$m=3$]{\includegraphics[scale=0.41]{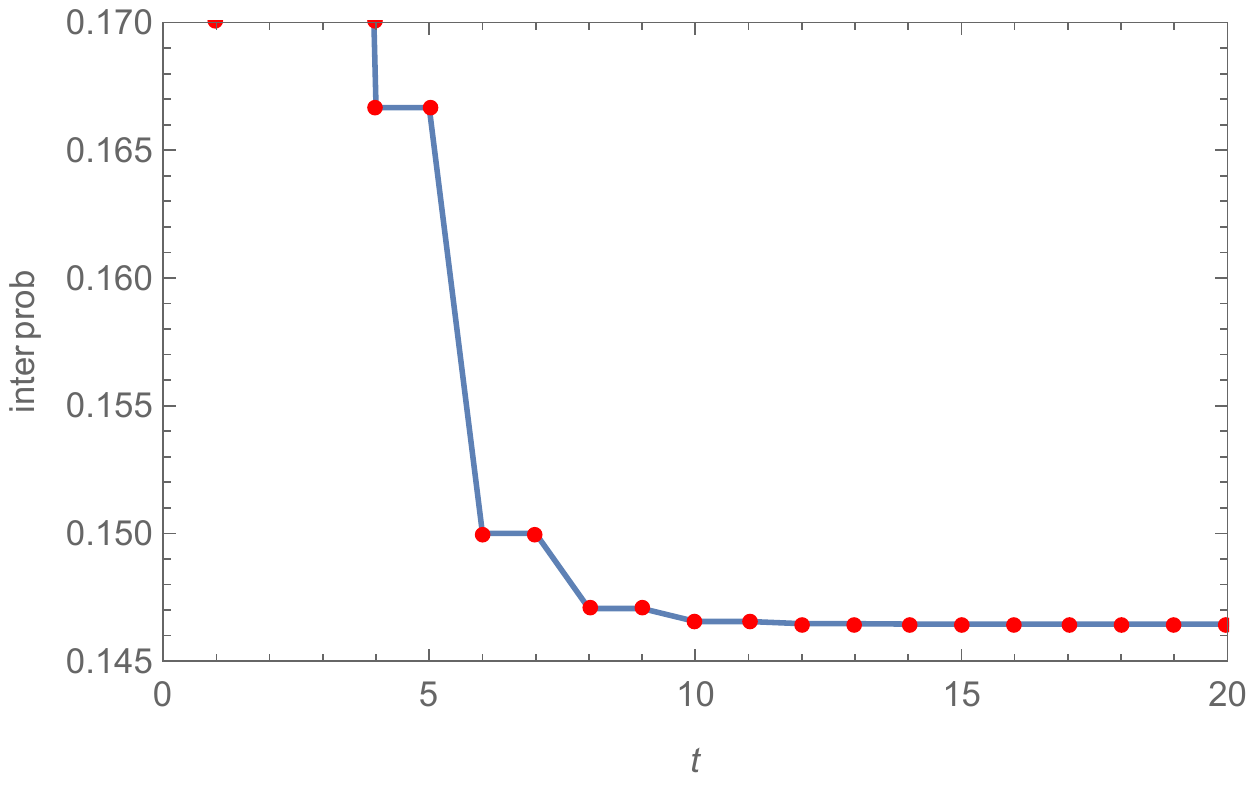}}
\hspace{1.75em}
\subfigure[$m=5$]{\includegraphics[scale=0.41]{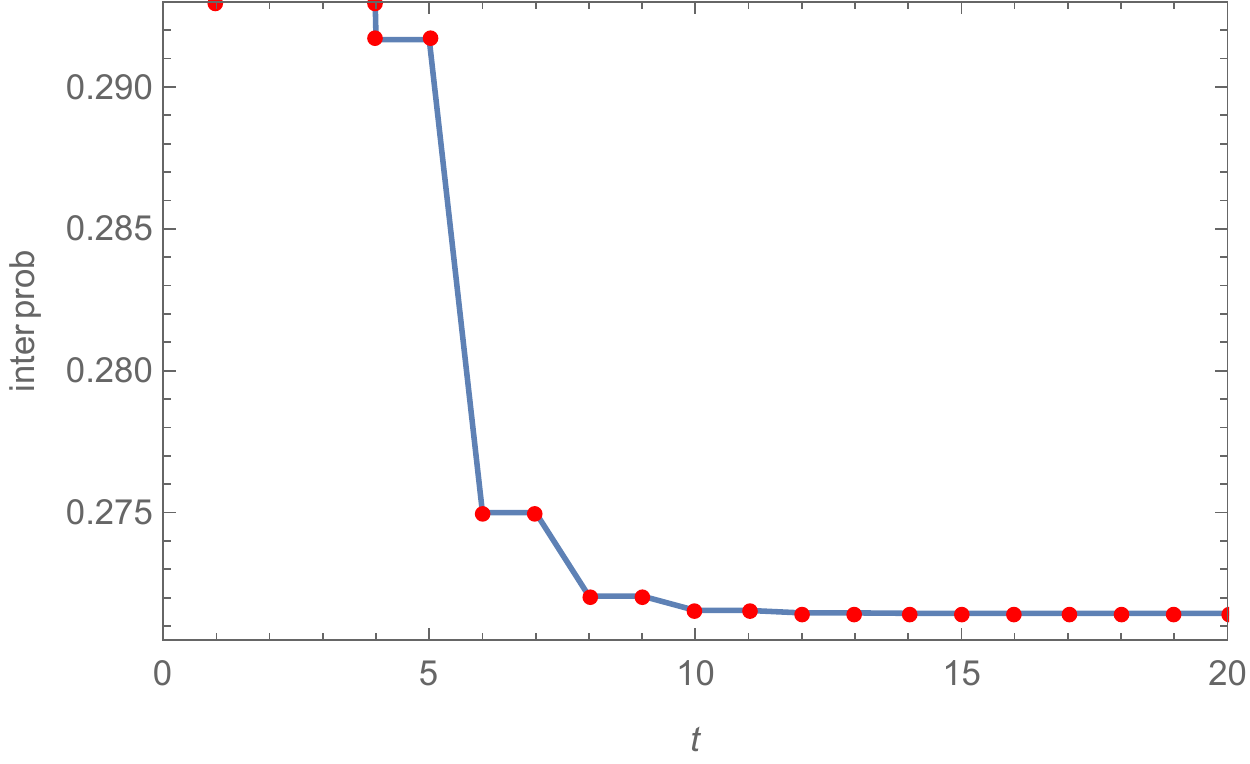}}
\caption{Interception probabilities for an attack at node $\mathit{1}$ of $L_{5}$, under the optimal patrol, for increasing delay.}
\label{Figure12}
\end{figure}

However, the same issues arise here regarding the validity of our findings for delay $d>15$, and the choice of the optimal attack node. To treat both we work the same way we did before. Firstly, we generate the values of $\pi_{m}(\hat{p},\hat{q},\hat{c},\hat{\kappa},d)$, $m=2,3,4,5,6$, for increasing values of $d$, and check if we reach the limiting interception probabilities without crossing below the optimal interception probabilities  we have estimated for $D=15$. As we see in Figure \ref{Figure12}, this appears to be the case for every $m$, which validates our findings, given in Table \ref{Table5}, also for $D>15$.

\begin{figure}[H]
\centering
\subfigure[Attack Node $\mathit{2}$]{\includegraphics[scale=0.4]{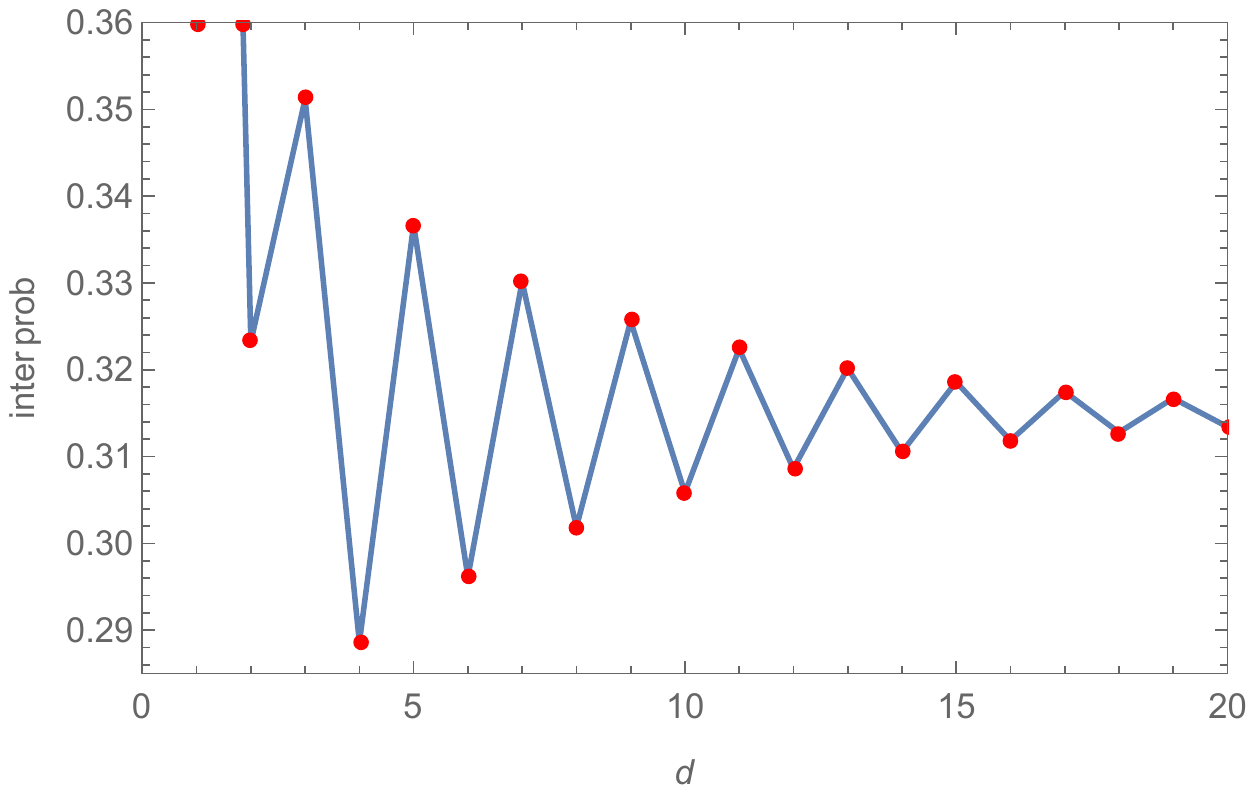}}
\hspace{1.75em}
\subfigure[Attack Node $\mathit{3}$]{\includegraphics[scale=0.4]{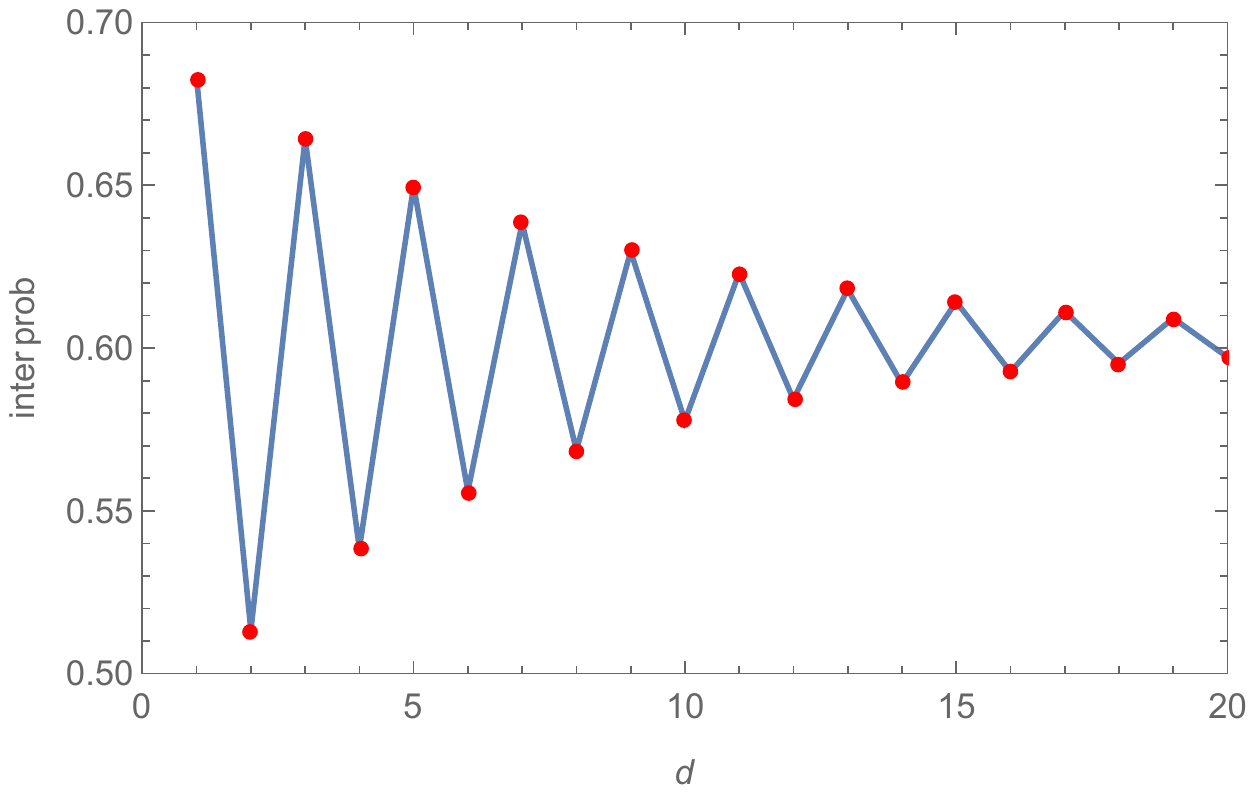}}
\caption{Interception probability for an attack at nodes $\mathit{2}$ and $\mathit{3}$ of $L_{5}$, under the optimal patrol estimated for an attack at node $\mathit{1}$, for increasing delay, for $m=4$.}
\label{Figure13}
\end{figure}

Regarding the choice of the optimal attack node, assuming the optimal patrols that we have estimated for an attack at node $\mathit{1}$, first we generate the corresponding interception probabilities for an attack at nodes $\mathit{2}$ and $\mathit{3}$ for increasing values of $d$, then we identify their minima occurring at a specific value of $d$, and last we compare those minima with the maxmin interception probabilities $\pi_{m}(\hat{p},\hat{q},\hat{c},\hat{\kappa},\hat{d})$ given in Table \ref{Table5}. For example, for $m=4$ we obtain Figure~\ref{Figure13}.

As seen in Figures \ref{Figure13}a and \ref{Figure13}b, for $m=4$, the minimum interception probabilities $0.2886$ and $0.5132$ that we find respectively for the attacks $(\mathit{2},4)$ and $(\mathit{3},2)$, are both greater than the minimum interception probability $0.1890$ that we have found for the attack $(\mathit{1},8)$ (all estimated for the optimal patrol $(0.4663,0.4330,0.4216,1)$, see Table \ref{Table5}). Similar results, confirming that the Attacker should comparatively attack an end node (node $\mathit{1}$ in this instance), we obtain for the rest of the attack durations $m$ we have considered, as we present in the Table \ref{Table6}.

\setlength{\tabcolsep}{13pt}
\renewcommand{\arraystretch}{0.8}
\begin{table}[H]
\centering
\begin{tabular}{ cccc } 
\hline\hline
 for $(\hat{p},\hat{q},\hat{c},\hat{\kappa})$ & $\pi_{m}$ at node $1$, $(i,\hat{d})$ & $\pi_{m}$ at node $2$, $(i,\hat{d})$ & $\pi_{m}$ at node $3$, $(i,\hat{d})$   \\ \hline
 $m=2$ &  $0.0646$, $(\mathit{1},10)$ & $0.0876$, $(\mathit{2},4)$ & $0.1080$, $(\mathit{3},2)$  \\ 
 $m=3$ &  $0.1464$, $(\mathit{1},12)$ & $0.2500$, $(\mathit{2},2)$ & $0.5000$, $(\mathit{3},1)$  \\ 
 $m=4$ &  $0.1890$, $(\mathit{1},8)$ & $0.2886$, $(\mathit{2},4)$ & $0.5132$, $(\mathit{3},2)$  \\ 
 $m=5$ &  $0.2714$,$(\mathit{1},14)$ & $0.4375$, $(\mathit{2},2)$ & $0.7500$, $(\mathit{3},1)$  \\ 
 $m=6$ &  $0.3000$, $(\mathit{1},8)$ & $0.4545$, $(\mathit{2},4)$ & $0.7473$, $(\mathit{3},2)$\\ \hline\hline
\end{tabular}
\caption{Compare attacks at nodes $\mathit{1}$ and $\mathit{2,3}$ under the optimal patrol for node $\mathit{1}$}
\label{Table6}
\end{table}

\section{The Circle Network $C_{n}$}\label{Section5}

The circle network $C_{n}$ is a closed walk consisting of $n$ nodes, and equal number of edges.

\begin{figure}[H]
\centering
\includegraphics[scale=0.85]{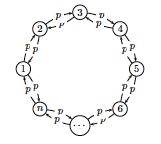}
\caption{The Circle network $C_{n}$.}
\label{Figure14}
\end{figure}

Like above, we restrict the Patroller to Markovian strategies such that he moves clockwise and counter clockwise with the same probability $p$, while he stays at the same node with probability $r=1-2\cdot p$. Given this symmetric patrol, every strategy $(i,d)$, for a random node $i$ and a given delay $d$, is equally advantageous for the Attacker. We note that an Associate Editor has observed that it is precisely our symmetry assumption (that the Patroller's motion respects any symmetries of the network) that prevents the Patroller from simply cycling clockwise to always intercept the attack when $m>n-1$. See Section \ref{Section2.1} for related comments.

According to our setting, the transition matrix $B_{n}$ characterizing the Patroller's Markovian walk on the circle network is parametrized by the single parameter $p$. Without loss of generality, we take node $\mathit{1}$ as the attack node, and we denote the \textit{Patroller's away distribution from node $\mathit{1}$} for $t\geq 1$ by $x^{(t)}=(0,x_{2}^{(t)},\dots,x_{n}^{(t)})$, where $x^{(0)}=\{1,0,\dots,0\}$.

Equivalently to our approach in Section \ref{Section4}, if we define the row vector for $t\geq 0$

\begin{equation*}
y^{(t)}=x^{(t)} \times B_{n},
\end{equation*}

\noindent then the Patroller's away distribution from node $\mathit{1}$ at $t\geq 1$ is given by the iteration formula

\begin{equation}\label{eq26}
x^{(t)}=\frac{x^{(t-1)}\times B_{n}-(y_{1}^{(t-1)}0,\dots,0)}{1-y_{1}^{(t-1)}}.
\end{equation}

We consider networks $C_{4}$ and $C_{5}$. We solve the game numerically for given values of $m$. Note that $C_{3}$ is the same as the complete graph $K_{3}$, and its solution is given by formula \eqref{eq34} given later in Section \ref{Section6.1}.

\subsection{The Network $C_{4}$}\label{Section5.1}
\indent

We start with the circle network with four nodes. Equivalently to \eqref{eq26}, we can rewrite the Patroller's away distribution for $C_{4}$ into the following system of equations

\begin{gather}\label{eq27}
\begin{aligned}
(1-p+p\cdot x_{3}^{(t-1)})\cdot x_{2}^{(t)}&=(1-2\cdot p)\cdot x_{2}^{(t-1)}+p\cdot x_{3}^{(t-1)},\\
(1-p+p\cdot x_{3}^{(t-1)})\cdot x_{3}^{(t)}&=p+(1-3\cdot p)\cdot x_{3}^{(t-1)}.
\end{aligned}
\end{gather}

System \eqref{eq27} defines a continuously differentiable map by $g_{2}(x_{2}^{(t)},x_{3}^{(t)})=(x_{2}^{(t+1)},x_{3}^{(t+1)})$. There is a unique stationary solution of $g_{2}$ of the form $(0,x_{2},x_{3},x_{4})$, independent of $p$, given~by

\begin{equation*}
(0,x_{2},x_{3},x_{4})=(0,\frac{2-\sqrt{2}}{2},\sqrt{2}-1,\frac{2-\sqrt{2}}{2}).
\end{equation*}

We consider four cases regarding the duration of the attack, $m=2,3,4,5$. The interception probability $\pi_{m}(p,d)$ in each case in given by

\begin{align*}
\pi_{2}(p,d)&=p\cdot (x_{2}^{(d)} + x_{4}^{(d)}),\\
\pi_{3}(p,d)&=2\cdot p\cdot (1-p)\cdot (x_{2}^{(d)}+x_{4}^{(d)}) + 2\cdot p^{2}\cdot x_{3}^{(d)},\\
\pi_{4}(p,d)&=p\cdot \bigl(2-2\cdot p+(1-2\cdot p)^{2}+2\cdot p^{2}\bigr)\cdot (x_{2}^{(d)}+x_{4}^{(d)}) + 2\cdot p^{2}\cdot (3-4\cdot p)\cdot x_{3}^{(d)},\\
\begin{split}
\pi_{5}(p,d)&=p\cdot \bigl(2-2\cdot p+(1-2\cdot p)^{2}+(1-2\cdot p)^{3}+2\cdot p^{2}\cdot (4-6\cdot p)\bigr)\cdot (x_{2}^{(d)}+x_{4}^{(d)})\\ &+2\cdot p^{2}\cdot \bigl(2\cdot p^{2}+2\cdot (1-2\cdot p)^{2}+3-4\cdot p\bigr)\cdot x_{3}^{(d)},
\end{split}
\end{align*}

We have estimated numerically, for $D=15$, the optimal game values for the above four cases, and we present our results in the following table rounded to four decimal places.

\setlength{\tabcolsep}{13pt}
\renewcommand{\arraystretch}{0.8}
\begin{table}[H]
\centering
\begin{tabular}{ ccccc } 
\hline\hline
 for $D=15$    & $\hat{p}$ & $\hat{d}$ & $\pi_{m}(\hat{p},\hat{d})$ & $\pi_{m}(\hat{p},\infty)$   \\ \hline
 $m=2$ & $0.2929$ & $2$ & $0.1716$ & $0.1716$ \\ 
 $m=3$ & $0.5$ & $1$ & $0.5$ & $0.5$  \\ 
 $m=4$ & $0.4515$ & $2$ & $0.5216$ & $0.6021$  \\
 $m=5$ & $0.5$ & $1$ & $0.75$ & $0.75$ \\ \hline\hline
\end{tabular}
\caption{Optimal game values for $C_{4}$}
\label{Table7}
\end{table}

\begin{figure}[H]
	\centering
\includegraphics[scale=0.41]{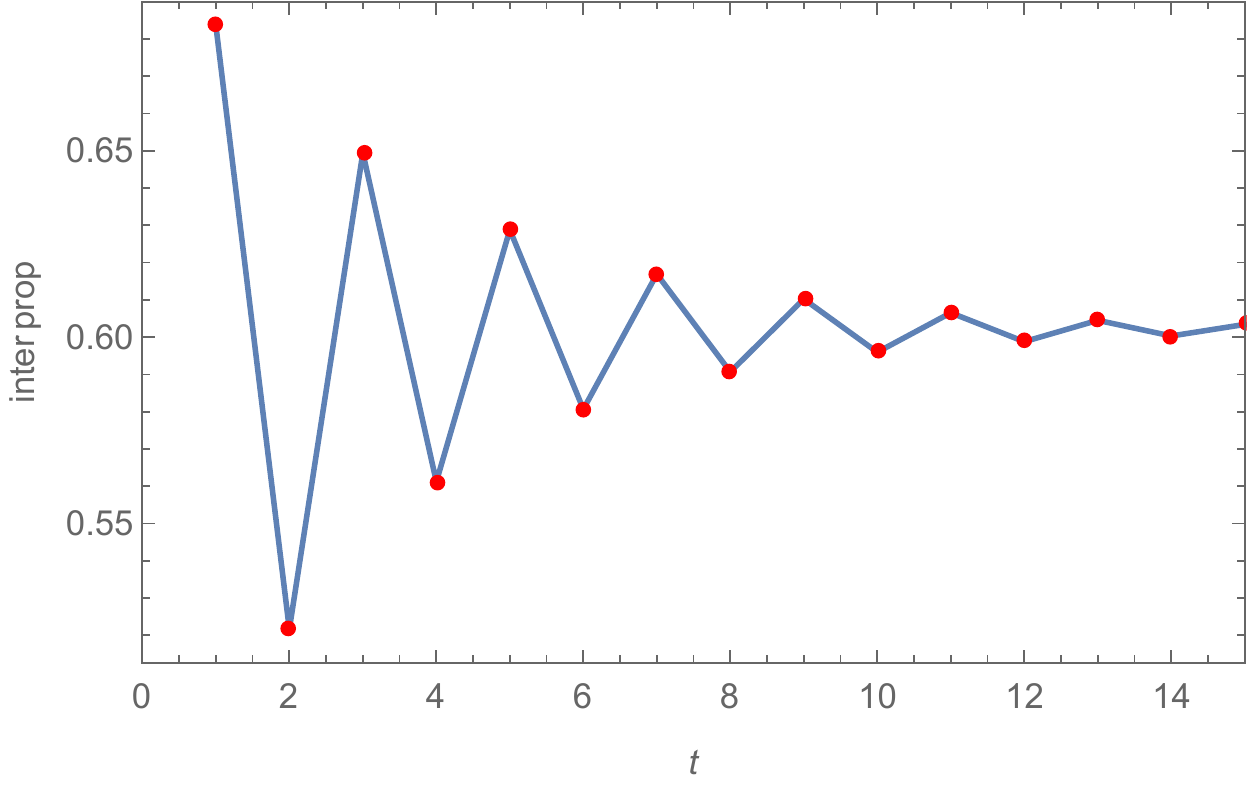}
	\caption{Interception probability for an attack at a random node of $C_{4}$ under the optimal patrol, for increasing delay, for $m=4$.}
	\label{Figure15}
\end{figure}

Again we want to examine (qualitatively) whether our numerical results are valid for $D>15$. For $m=2,3,5$, the optimal interception probabilities coincide with the corresponding limiting interception probabilities. Similarly for $m=4$, we reach the limiting interception probability without crossing below the optimum interception probability that we found for $D=15$, as we see in Figure \ref{Figure15}. Thus, we may claim that our results are valid also for infinitely large delay $d$.

\subsection{The Network $C_{5}$}\label{Section5.2}
\indent

Equivalently to \eqref{eq26}, we rewrite the Patroller's away distribution from node $\textit{1}$ into the following system of equations

\begin{gather}\label{eq29}
\begin{aligned}
(1-p+p\cdot x_{3}^{(t-1)}+ p\cdot x_{4}^{(t-1)})\cdot x_{2}^{(t)} &= (1-2\cdot p)\cdot x_{2}^{(t-1)}+p\cdot x_{3}^{(t-1)},\\
(1-p+p\cdot x_{3}^{(t-1)}+ p\cdot x_{4}^{(t-1)})\cdot x_{3}^{(t)} &= p\cdot (x_{2}^{(t-1)}+x_{4}^{(t-1)})+(1-2\cdot p)\cdot x_{3}^{(t-1)},\\
(1-p+p\cdot x_{3}^{(t-1)}+ p\cdot x_{4}^{(t-1)})\cdot x_{4}^{(t)} &=p+(1-3\cdot p)\cdot x_{4}^{(t-1)}-p\cdot x_{2}^{(t-1)}.
\end{aligned}
\end{gather}

System \eqref{eq29} defines a continuously differentiable mapping $g_{3}:\Delta^{3}\mapsto\Delta^{3}$. There is a unique stationary solution of $g_{3}$ of the form $(0,x_{2},x_{3},x_{4},x_{5})$, that is independent of $p$, given by 

\begin{equation*}
(0,x_{2},x_{3},x_{4},x_{5})=(0,\frac{3-\sqrt{5}}{4},\frac{\sqrt{5}-1}{4},\frac{\sqrt{5}-1}{4},\frac{3-\sqrt{5}}{4}).
\end{equation*}

We consider here the same four cases regarding the duration $m$ of the attack, and we have respectively the following interception probabilities

\begin{align*}
\pi_{2}(p,d)&=p\cdot (x_{2}^{(d)}+x_{5}^{(d)}),\\
\pi_{3}(p,d)&=2\cdot p\cdot (1-p)\cdot (x_{2}^{(d)}+x_{5}^{(d)})+p^{2}\cdot (x_{3}^{(d)}+x_{4}^{(d)}),\\
\pi_{4}(p,d)&=p\cdot \bigl(2-2\cdot p+(1-2\cdot p)^{2}+p^{2}\bigr)\cdot (x_{2}^{(d)}+x_{5}^{(d)})
+3\cdot p^{2}\cdot (1-p)\cdot (x_{3}^{(d)}+x_{4}^{(d)}),\\
\begin{split}
\pi_{5}(p,d)&=p\cdot \bigl(2-2\cdot p +(1-2\cdot p)^{2}+p^{2}+2\cdot p^{2}\cdot (1-2\cdot p)+p^{3}\bigr)\cdot (x_{2}^{(d)}+x_{5}^{(d)})\\
&+p^{2}\cdot \bigl(3-4\cdot p^{2}+3\cdot (1-2\cdot p)^{2}\bigr)\cdot (x_{3}^{(d)}+x_{4}^{(d)}).
\end{split}
\end{align*}

We have estimated numerically, for $D=15$, the optimal game values for the above four cases, and we gather our results in the following table rounded to four decimal places.

\setlength{\tabcolsep}{13pt}
\renewcommand{\arraystretch}{0.8}
\begin{table}[H]
\centering
\begin{tabular}{ ccccc } 
\hline\hline
 for $D=15$    & $\hat{p}$ & $\hat{d}$ & $\pi_{m}(\hat{p},\hat{d})$ & $\pi_{m}(\hat{p},\infty)$   \\ \hline
 $m=2$ & $0.3820$ & $2$ & $0.1459$  & $0.1459$ \\ 
 $m=3$ & $0.3820$ & $2$ & $0.2705$  & $0.2705$  \\ 
 $m=4$ & $0.4450$ & $2$ & $0.3808$  & $0.4282$  \\ 
 $m=5$ & $0.5$ & $2$ & $0.5$  & $0.5716$ \\ \hline\hline
\end{tabular}
\caption{Optimal game values for $C_{5}$}
\label{Table8}
\end{table}

We want to check whether setting $D>15$ (for $p=\hat{p}$), we reach the limiting interception probability without crossing below the optimum interception probability estimated for $D=15$. 

\begin{figure}[H]
	\centering
	\subfigure[$m=4$]{\includegraphics[scale=0.41]{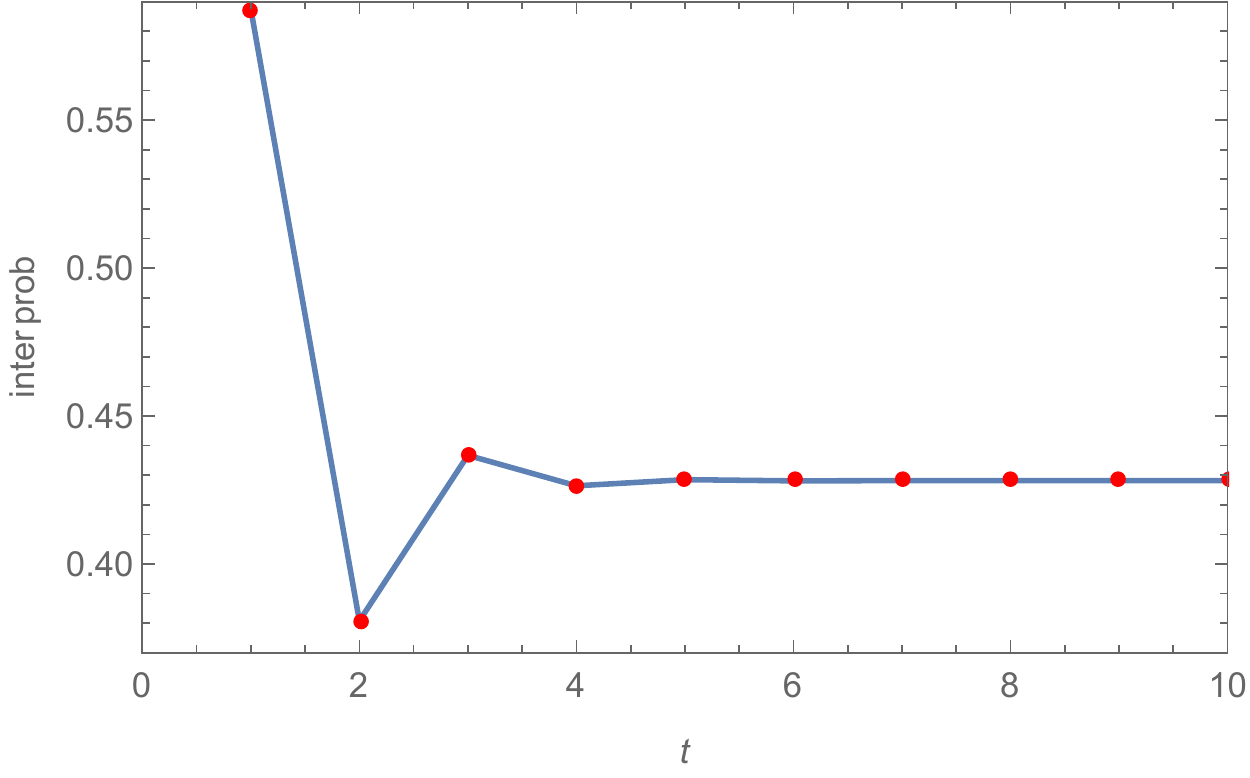}}\hspace{1.75em}
           \subfigure[$m=5$]{\includegraphics[scale=0.41]{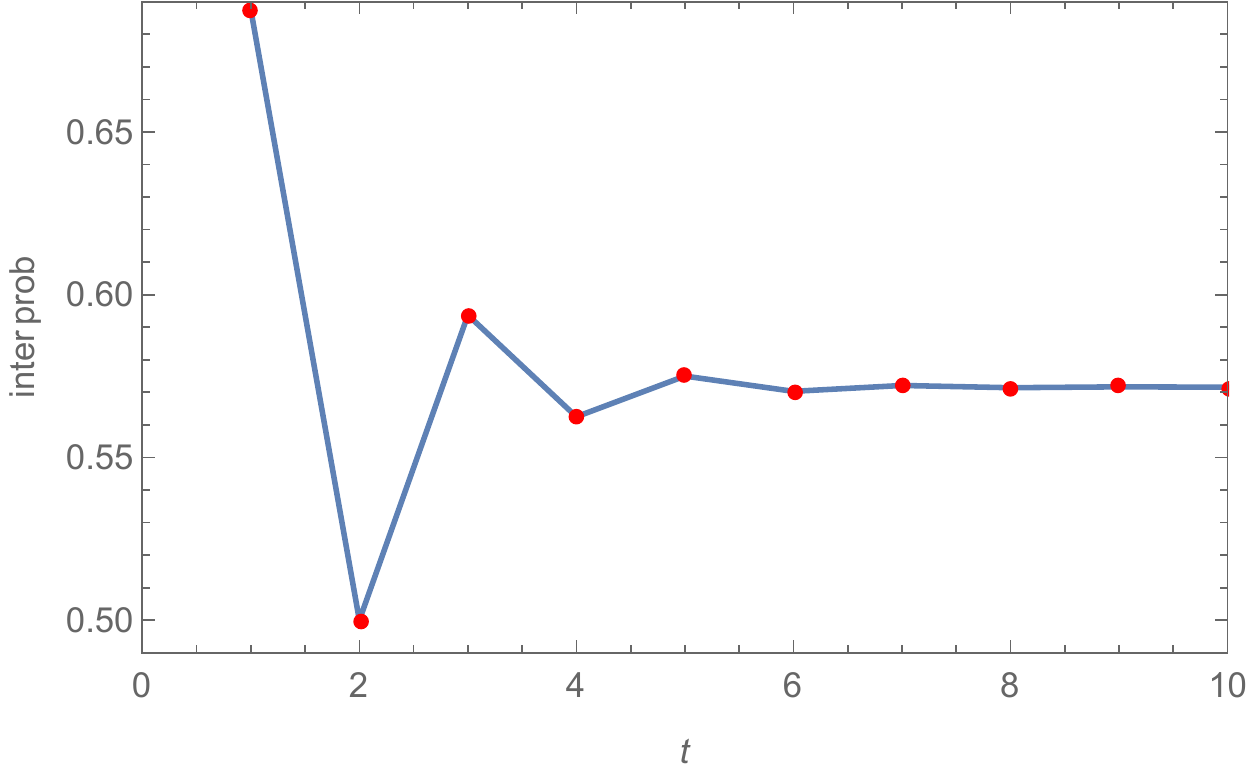}}
	\caption{Interception probability for an attack at a random node of $C_{5}$ under the optimal patrol, for increasing delay.}
	\label{Figure16}
\end{figure}

\noindent For $m=4,5$ this appears to be the case as seen in Figure \ref{Figure13}, while for $m=2,3$ the optimum interception probability coincides with the limiting interception probability.

\section{The Star-in-Circle Network $E_{n}$}\label{Section6}
\indent

The Star-in-Circle network $E_{n}$ is a combination of the star $S_{n}$ and the circle $C_{n}$ networks we have considered above. It consists of $n$ end nodes, a single central node, and $2\cdot n$ edges.

\begin{figure}[H]
\centering
\includegraphics[scale=0.85]{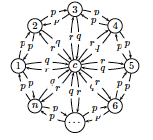}
\caption{The Star-in-Circle network $S_{n}$.}
\label{Figure17}
\end{figure}

We take a Markovian Patroller that leaves an end towards each of its adjacent ends with probability $p$ and towards the center with probability $q$, leaves the center towards each end with probability $r$, remains at an end with probability $a=1-2\cdot p-q$, and at the center with probability $b=1-n\cdot r$. The transition matrix $\Gamma_{n}$ characterizing his walk on $E_{n}$ is parametrized by  $p$, $q$, $r$. Regarding the Attacker, she can either attack at one of the ends, that are equally preferable, or at the center. 

Without loss of generality, we take node $\textit{1}$ as the end node to be attacked. In accordance with \eqref{eq18} and \eqref{eq26}, the Patroller's away distribution from node $\textit{1}$ for $t\geq 1$ is given by the formula

\begin{equation}\label{eq31}
x^{(t)}=\frac{x^{(t-1)}\times \Gamma_{n}-(z_{1}^{(t-1)},0,\dots,0)}{1-z_{1}^{(t-1)}},
\end{equation}

\noindent where $z^{(t)}=x^{(t)}\times \Gamma_{n}$ is the row vector defined for $t\geq 0$, and $x^{(0)}=(1,0,\dots,0)$.

The Patroller's away distribution from the center for $t\geq 1$ is defined by $y^{(t)}=(y_{1}^{(t)},\dots,y_{n}^{(t)},0)$. If we consider an attack at the center of duration $m$, then the interception probability is given by

\begin{equation}\label{eq32}
\pi_{c}(q,m)=\sum_{n=2}^{m}q\cdot (1-q)^{n-2}.
\end{equation}

Recall that our aim is to find the optimal Patroller's walk $(\hat{p},\hat{q},\hat{r})$, and the optimal Attacker's waiting time $\hat{d}$. First though we need to examine which node the Attacker will optimally attack. Obviously, if the Patroller announces a $\hat{q}$ such that $\pi_{c}(\hat{q},m)<\pi_{m}(p,\hat{q},r,d)$ then the Attacker should attack the center. Therefore, here we deal with the following problem

\begin{equation}\label{eq33}
\max_{\Gamma_{n}}\min_{d}\pi_{m}(p,q,r,d)=\pi_{m}(\hat{p},\hat{q},\hat{r},\hat{d}) : \pi_{m}(\hat{p},\hat{q},\hat{r},\hat{d})\leq \pi_{c}(\hat{q},m),
\end{equation}

\noindent the solution of which actually makes the Attacker indifferent of which node to attack.

\subsection{The Network $E_{3}$ and The Complete Network $K_{n}$}\label{Section6.1}
\indent

First we consider the star-in-circle network $E_{3}$ with three end nodes. This is essentially the same as the complete network $K_{4}$, except that in $E_{3}$ there is a distinguished node, that is the center $c$. The network $K_{n}$ is relatively simple to analyze in for our game. At the first period of an attack of duration $m$, the Patroller is away from the attack node $e$. In each subsequent period, assuming a random walk, the probability that he goes to a different node (other than $e$) is $(n-2)/(n-1)$, so the probability of intercepting the attack is given by 

\begin{equation}\label{eq34}
1-((n-2)/(n-1))^{m-1}.
\end{equation}

\noindent Note that for $n=4$ and $m=2,3,4$, this gives $1/3$, $5/9=0.5556$, and $19/27=0.7037$ respectively.

In fact, the optimal interception probability for the complete graph $K_{n}$ with duration $m$ is given by formula \eqref{eq34} because, if the Patroller stays at a node with a positive probability, the interception probability is strictly less than that given by \eqref{eq34}. The Attacker can choose any delay $d$.

\subsection{The Network $E_{4}$}\label{Section6.2}
\indent

Next we consider the star-in-circle network with four end nodes. The interception probability of an attack of duration $m=2,3,4$, and delay $d$, at node $\textit{1}$ is given respectively by

\begin{align*}
\pi_{2}(p,q,r,d)&=p\cdot (x_{2}^{(d)}+x_{4}^{(d)})+r\cdot x_{c}^{(d)},\\
\pi_{3}(p,q,r,d)&=(p+a\cdot p+q\cdot r)\cdot (x_{2}^{(d)}+x_{4}^{(d)})
+(2\cdot p^{2}+q\cdot r)\cdot x_{3}^{(d)}+r\cdot (1+2\cdot p+b)\cdot x_{c}^{(d)},\\
\begin{split}
\pi_{4}(p,q,r,d)&=\bigl(2\cdot p^{3}+3\cdot q\cdot r\cdot p+p+a\cdot r+q\cdot r+a^{2}\cdot p+a\cdot q\cdot r+q\cdot b\cdot r\bigr)\cdot(x_{2}^{(d)}+x_{4}^{(d)})\\&+\bigl(2\cdot p^{2}+q\cdot r+4\cdot a\cdot p^{2}+4\cdot p\cdot q\cdot r+a\cdot q\cdot r+q\cdot b\cdot r\bigr)\cdot x_{3}^{(d)}\\&+r\cdot \bigl(1+2\cdot p+b+2\cdot b\cdot p+2\cdot a\cdot p+b^{2}+3\cdot r\cdot q+2\cdot p^{2}\bigr)\cdot x_{c}^{(d)}.
\end{split}
\end{align*}

We have solved the game numerically, for $D=15$, and we gather our results in the following table rounded to four decimal places.

\setlength{\tabcolsep}{13pt}
\renewcommand{\arraystretch}{0.8}
\begin{table}[H]
\centering
\begin{tabular}{ cccccc } 
\hline\hline
 for $D=15$    & $(\hat{p},\hat{q},\hat{r})$ & $\hat{d}$ & $\hat{\pi}_{m}(\hat{p},\hat{q},\hat{r},\hat{d})$    \\ \hline
 $m=2$ & $(0.2835,0.1695,0.25)$ & $2$ & $0.1695$ \\ 
 $m=3$ & $(0.3886,0.2229,0.25)$ & $2$ & $0.3961$  \\ 
 $m=4$ & $(0.3945,0.2109,0.25)$ & $2$ & $0.5087$  \\ \hline\hline
\end{tabular}
\caption{Optimal game values for $E_{4}$}
\label{Table9}
\end{table}

\subsection{Giving the Attacker greater vision}\label{Section6.3}
\indent

Our basic model restricts the Attacker's vision only to observing when the Patroller is present at her attack node and when he is not. Here, at the suggestion of an anonymous referee, we extend the Attacker's vision, and we examine for example if the Attacker can benefit by observing the edge by which the Patroller leaves her attack node. To keep things as simple as possible otherwise, we do this in the setting of the star-in-circle network $E_{4}$, with difficulty $m=2$. Note that since the case of $E_3$ is essentially the complete graph $K_{4}$, there is too much symmetry to distinguish the leaving edge.

From the point of view of the Attacker at the end node $e$, we divide the states, namely the Patroller's possible locations, into three types: $C$ (for the center), $A$ (for an adjacent end node), and $O$ (for the opposite end node). Clearly the Attacker would prefer to initiate the attack when the Patroller is in state $O$, as this would result in a certain successful attack for $m=2$. The transition probabilities, given that the Patroller stays away from the Attacker's node $e$, are given by the matrix $M$ (see Figure \ref{Figure17})

\begin{equation}\label{eq36}
\bordermatrix{
\mbox{states} & C & A & O \cr
C & \frac{b}{1-r} & \frac{2r}{1-r} & \frac{r}{1-r} \cr
A & \frac{q}{1-p} & \frac{a}{1-p} & \frac{p}{1-p} \cr
O & q & 2p & a \cr
}\equiv M.
\end{equation}

\noindent If the Patroller leaves the attack node $e$ to go to an adjacent end node, then the sequence of distributions over $C$, $A$, and $O$ are given by

\begin{equation*}
x_{A}^{t}=(0,1,0)\cdot M^{t},
\end{equation*}

\noindent while, if the Patroller leaves $e$ to go to the center, then the sequence of distributions are given by

\begin{equation*}
x_{C}^{t}=(1,0,0)\cdot M^{t}.
\end{equation*}

According to what is the state of the Patroller when the attack at $e$ starts, the interception probability is $r$, $q$, or $0$, for the states $C$, $A$, and $O$ respectively. Assuming the distribution over the three states is $x=(x_{1},x_{2},x_{3})$ when the attack starts, the interception probability is given by $(r,q,0)\cdot x=r\cdot x_{1}+q\cdot x_2$. Moreover, if the Attacker starts her attack at an end $d_{A}$ periods after the Patroller goes to an adjacent end or $d_{C}$ periods after the Patroller goes to the center, or if she attacks at the center at any time (given the Patroller is not there), then the three interception probabilities are given respectively by

\begin{equation*}
(r,q,0)\cdot x_{A}^{d}\quad,\quad (r,q,0)\cdot x_{C}^{d},\quad \text{and}\quad q.
\end{equation*}

\noindent Hence, the problem for the Patroller has the following form

\begin{equation}\label{eq37}
\max_{p,q,r}\min_{d_{A},d_{C}\leq D}\bigl((r,q,0)\cdot x_{A}^{d},(r,q,0)\cdot x_{C}^{d}\bigr).
\end{equation}

Solving \eqref{eq37} numerically we find that the interception probability is approximately $0.1667$ for what is the Patroller's optimal walk $(p,q,r)=(0.25,0.1667,0.25)$, with $\hat{d}_{A}=\hat{d}_{C}=2$. Thus, by letting the Attacker observe the edge by which the Patroller leaves the attacked node, the interception probability is reduced from the original value of $0.1695$ (see Table \ref{Table9}) to the new value of $0.1667$, a reduction of approximately $(0.1695-0.1667)/0.1695\simeq 1.7\%$. Against the previously optimal Patroller's strategy of $(p,q,r)=(0.2835,0.1695,0.25)$ the respective interception probabilities when the attack is after the Patroller leaves by an adjacent node or the central node are about about $0.18$ and $0.13$. So the Attacker's optimal response would be to attack when the Patroller leaves by the central node.

To obtain a clearer view we implement more algebraic methods to investigate \eqref{eq37}. Based on our numerical result of $\hat{d}_{A}=\hat{d}_{C}=2$, the interception probabilities of an attack at an end after the Patroller leaves via the center or an adjacent end are respectively

\begin{equation}\label{eq38}
\frac{r\cdot (1+2\cdot p-4\cdot r)}{1-r},
\end{equation} 
\begin{equation}\label{eq39}
\frac{p\cdot (1-q)+q\cdot r-2\cdot p^2}{1-p}.
\end{equation}

\noindent Equating \eqref{eq38} with \eqref{eq39} and solving in terms of $p$ we get $p=r$. Substituting then back to \eqref{eq38} we come up with the interception probability in the form

\begin{equation*}
f(r)=\frac{r\cdot(1-2\cdot r)}{1-r},
\end{equation*}

\noindent which is maximized when 

\begin{equation*}
f^{'}(r)=\frac{1-4r+2r^2}{(1-r)^2}=0,\quad \text{or} \quad r=\hat{r}=\hat{p}=1/4,
\end{equation*}

\noindent with $f(1/4)=1/6\simeq 0.1667$ as the interception probability.

If the Patroller chooses any $q<1/6$, then since $q$ is the interception probability for any attack at the center, the Attacker would respond by attacking at the center. However, the Patroller can ensure an interception probability of at least $1/6$ by choosing any value $1/6\leq q\leq 1/2$. For $q=1/6$ all three types have the same interception probability.

\section{Conclusion}\label{Section7}
\indent 

This paper introduces limited vision to the Attacker in a patrolling game by letting her observe, at the beginning of every period, whether or not the Patroller is present at the node she is planning to attack. Of course the Attacker will not begin the attack when the Patroller is present, and we show here that optimally she will also not immediately attack when the Patroller leaves but rather waits an optimal number of periods, resetting her count if the Patroller returns before that time. This behavior is in stark contrast to the optimal behavior when the Patroller follows a known periodic tour of the network, where an immediate attack is often the best strategy. For example many prison escape movies show the prisoners (i.e. the Attackers) attempting to escape just after the spotlight (i.e. the Patroller) leaves their location.

We have adopted the assumption that the Patroller's motion can be observed prior to the game, namely that is somehow a known element of the game. That is, we take a Stackelberg approach. However, an interesting observation of our results for several networks is that in fact the Attacker has an optimal strategy that does not require any knowledge of the Patroller's motion, namely a Nash equilibrium exists between the Patroller and the Attacker.

A number of possible extensions of our model naturally suggest themselves after our results are examined. For one, it would be interesting to give the Attacker a greater range of vision, perhaps initially to nodes adjacent to the planned attack node. Of course, if the Attacker has significant vision in this respect, he might have an incentive to move away from the Patroller's current location by changing his choice of which node to attack. This approach might lead to one-sided information in the so called Cops and Robbers games, with only the Robbers (our Attacker) having vision. In terms of search and pursuit-evasion games wording, this might be called a search-evasion game approach. Another extension to our model, seemingly rather difficult, would be to allow the Patroller the same range of motions (mixtures over general walks) that was allowed in the original formulation of patrolling games \cite{Alp-1}. A first step towards this direction would be to introduce Markovian strategies with short memories.

Moreover, our model raises the question of why patrols are ever carried out with uniforms that give a potential attacker or infiltrator additional helpful information. One answer might lie in the direction of deterrence. Under some of our game parameters, where the interception probability is relatively high, the Attacker might choose to abandon his attempt altogether, leading to a non zero-sum game. Clearly, our model of a Uniformed Patroller leads to many unanswered questions.

\section*{Acknowledgements}
Steve Alpern wishes to acknowledge support from the US Air Force Office of Scientific Research under grant FA9550-14-1-0049. Stamatios Katsikas wishes to acknowledge support from the Engineering and Physical Sciences Research Council (EPSRC), Swindon, UK.

\section*{Conflicts of Interest}
The authors declare no conflict of interest.

\end{document}